# SHORTEST SPANNING TREES AND A COUNTEREXAMPLE FOR RANDOM WALKS IN RANDOM ENVIRONMENTS

By Maury Bramson,[1] Ofer Zeitouni[2] and Martin P. W. Zerner[3]

*University of Minnesota, Technion and Universität Tübingen*

We construct forests that span $\mathbb{Z}^d$, $d \geq 2$, that are stationary and directed, and whose trees are infinite, but for which the subtrees attached to each vertex are as short as possible. For $d \geq 3$, two independent copies of such forests, pointing in opposite directions, can be pruned so as to become disjoint. From this, we construct in $d \geq 3$ a stationary, polynomially mixing and uniformly elliptic environment of nearest-neighbor transition probabilities on $\mathbb{Z}^d$, for which the corresponding random walk disobeys a certain zero–one law for directional transience.

**1. Introduction.** Let $d \geq 2$ and let $a : \mathbb{Z}^d \to \mathbb{Z}^d$ be a random function for which $x$ and $a(x)$ are always nearest neighbors. If $a(a(x)) \neq x$ for all $x$ and the set $F_a = \{\{x, a(x)\} | x \in \mathbb{Z}^d\}$ of edges defines a forest in $\mathbb{Z}^d$ [i.e., the graph $(\mathbb{Z}^d, F_a)$ does not have cycles], we call such a random function $a$ an *ancestral function*. In particular, if $a$ is an ancestral function, then each connected component of $F_a$ is infinite and we can interpret $a(x)$ as the parent or immediate ancestor of $x$. The $n$th generation ancestor of $x$, $n \geq 1$, is denoted by $a^n(x) = a(a^{n-1}(x))$, where $a^0(x) = x$. Then, Ray$(x) = \{a^n(x) | n \geq 0\}$ is the set of ancestors of $x$, including $x$ itself, whereas Tree$(x) = \{y \in \mathbb{Z}^d | x = a^n(y) \text{ for some } n \geq 1\}$ is the set of progeny of $x$. The length of the longest branch in Tree$(x)$ is defined as

(1)    $h(x) = \sup\{n \geq 0 | x = a^n(y) \text{ for some } y \in \mathbb{Z}^d\}, \qquad x \in \mathbb{Z}^d$.

Received January 2005; revised June 2005.
[1] Supported in part by NSF Grant DMS-02-26245.
[2] Supported in part by NSF Grant DMS-03-02230.
[3] Most of this paper was written while M. Z. was at Stanford University and at University of California, Davis.

In this paper, we study the tail behavior of $h(0)$ for such forests $F_a$ that are also stationary with respect to the translations of the lattice $\mathbb{Z}^d$ [that is, the collection $(a(x) - x)_{x \in \mathbb{Z}^d}$ is stationary]. The proof of the following theorem is easy.

THEOREM 1. *There is a constant $c_1 > 0$ that depends only on $d \geq 2$, such that for all stationary ancestral functions $(a(x))_{x \in \mathbb{Z}^d}$,*

$$\liminf_{n \to \infty} n^{d-1} \mathbb{P}[h(0) \geq n] \geq c_1. \tag{2}$$

Here and throughout the paper, $\mathbb{P}$ will denote the probability measure of the underlying probability space. The corresponding expectation operator will be denoted by $\mathbb{E}$.

An ancestral function $(a(x))_{x \in \mathbb{Z}^d}$ is *directed* if for some $z \in \{\pm 1\}^d$, $a(x) - x \in \{z_i e_i | i = 1, \ldots, d\}$ for all $x \in \mathbb{Z}^d$, $\mathbb{P}$-a.s., where $e_1, \ldots, e_d$ denote the standard basis vectors of $\mathbb{R}^d$. We then refer to $z$ as the direction of $a$ (or of the corresponding forest). Perhaps the simplest example of a stationary directed forest spanning $\mathbb{Z}^d$ is given in the following example.

EXAMPLE 1. Define $a(x) = x + i(x)$, where $i(x)$, $x \in \mathbb{Z}^d$, are independent random variables with $\mathbb{P}[i(x) = e_j] = 1/d$ for $j = 1, \ldots, d$. This defines a directed forest that spans $\mathbb{Z}^d$. Part of such a forest is shown in Figure 1(a), in $d = 2$. (It is not difficult to show that, in $d = 2$, the forest consists of a single tree, $\mathbb{P}$-a.s.; see [9].) It follows from the discussion in [9], page 1730, that Tree(0) is enclosed by two directed simple symmetric random walk paths on the dual lattice that are independent of each other until they meet. So, $\mathbb{P}[h(0) \geq n] \geq c_2 n^{-1/2}$ for some $c_2 > 0$ and all $n \geq 1$. (Neither this last fact nor Example 1 is used in the sequel, except as motivation.)

Example 1 might suggest that trees in stationary spanning forests need to be longer than suggested in Theorem 1, that is, that the rate of decay given in Theorem 1 is not optimal. However, this is not the case, as is shown by the following result.

THEOREM 2. *For each $d \geq 2$, there is a stationary and directed ancestral function $(a(x))_{x \in \mathbb{Z}^d}$ that is polynomially mixing of order 1 and for which*

$$\limsup_{n \to \infty} n^{d-1} \mathbb{P}[h(0) \geq n] < \infty. \tag{3}$$

Here, we are using the following notion of mixing.

DEFINITION 1. Let $b = (b(y))_{y \in \mathbb{Z}^d}$ be a family of random variables on some common probability space. For $G \subset \mathbb{Z}^d$, define the collections of real-valued random variables

(4) $\mathcal{M}_G^b = \{f : |f| \leq 1, f \text{ is measurable with respect to } \sigma(b(y), y \in G)\}.$

For a given $\gamma > 0$, $b$ is *polynomially mixing* (*of order $\gamma$*) if for all finite $G \subset \mathbb{Z}^d$,

$$\sup_{s \in \mathbb{Z}^d} \sup_{f \in \mathcal{M}_G^b, g \in \mathcal{M}_{G+s}^b} |s|^\gamma |\operatorname{cov}(f,g)| < \infty.$$

Our motivation for studying the above growth properties of random forests in $\mathbb{Z}^d$ was our desire to investigate possible extensions of a conjectured 0–1 law for random walk in a random environment (RWRE). We proceed to introduce the RWRE model.

For $d \geq 1$, let $S$ denote the set of $2d$-dimensional probability vectors and set $\Omega = S^{\mathbb{Z}^d}$. We consider all $\omega \in \Omega$, written as $\omega = ((\omega(x, x+e))_{|e|=1})_{x \in \mathbb{Z}^d}$, as an *environment* for the random walk that we define next. The random walk in the environment $\omega$, started at $z \in \mathbb{Z}^d$, is the Markov chain $(X_n)_{n \geq 0}$ with state space $\mathbb{Z}^d$, such that $X_0 \equiv z$, and whose transition probabilities $P_\omega^z$ satisfy

(5) $P_\omega^z(X_{n+1} = x + e | X_n = x) = \omega(x, x+e) \quad \text{for } e \in \mathbb{Z}^d \text{ with } |e| = 1.$

An environment $\omega$ is called *elliptic* if $\omega(x, x+e) > 0$ for all $x, e \in \mathbb{Z}^d$ with $|e| = 1$. A random environment $\omega$ is called *uniformly elliptic* if there exists a so-called *ellipticity constant* $\kappa > 0$, such that $\mathbb{P}[\omega(x, x+e) > \kappa] = 1$ for all $x, e \in \mathbb{Z}^d$ with $|e| = 1$. (See [6] and [8] for an introduction to the RWRE model and its properties.)

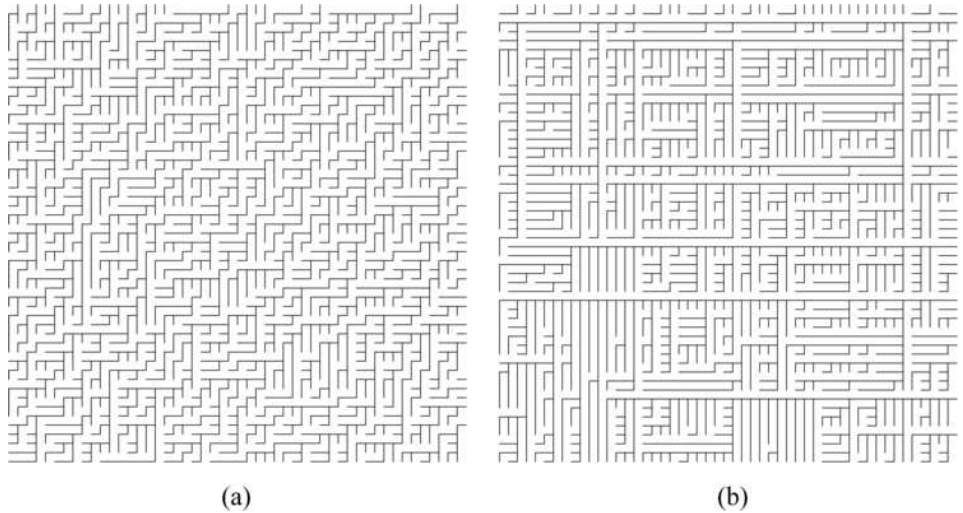

FIG. 1. (a) *Example* 1. (b) *Part of the forest constructed for $d = 2$ in the proof of Theorem* 2. *Note the long straight branches in the latter case.*

One of the major open questions in the study of the RWRE concerns the so-called 0–1 *law*. Fix a vector $\ell \in \mathbb{R}^d$, $\ell \neq 0$, and define the events $A_+(\ell)$ and $A_-(\ell)$ by

$$A_\pm(\ell) = \left\{ \lim_{n \to \infty} X_n \cdot \ell = \pm\infty \right\}.$$

It has been known since the work of Kalikow [2] that if the random vectors $\omega(x, \cdot)$, $x \in \mathbb{Z}^d$, are i.i.d. and $\omega$ is uniformly elliptic, then

$$\mathbb{E}[P_\omega^0[A_+(\ell) \cup A_-(\ell)]] \in \{0, 1\}.$$

This was extended in [9], Proposition 3, to the elliptic i.i.d. case.

The 0–1 law conjecture for RWRE states that if $\omega(x, \cdot)$, $x \in \mathbb{Z}^d$, are i.i.d. and $\omega$ is uniformly elliptic, then, in fact,

(6) $$\mathbb{E}[P_\omega^0[A_+(\ell)]] \in \{0, 1\}.$$

(The origin of this conjecture is a bit murky. For $d = 1$, it is a consequence of the law of large numbers in [5]. Kalikow [2] presented it as a question in $d = 2$; that case was settled only recently in the affirmative in [9]. The conjecture has since become folklore and is mentioned, e.g., in [7]. Although the question has arisen whether (6) holds for elliptic i.i.d. environments or for uniformly elliptic ergodic environments, we state the conjecture here in the weaker form, that is, for uniformly elliptic i.i.d. environments. For $d \geq 3$, this is still an open problem.) Recently, it has been shown that the law of large numbers for the RWRE follows from (6) for i.i.d. environments [10], and for a class of Gibbsian environments [3].

When $d = 2$ and the environment is elliptic and i.i.d., (6) was proved in [9], using techniques that do not extend to higher dimensions. The same paper provides an example, based on a construction of a forest that spans $\mathbb{Z}^2$, of an elliptic, ergodic environment, where (6) fails. However, this environment is neither uniformly elliptic nor mixing (in the ergodic theoretic sense), and not even totally ergodic, and thus the results in [9] do not contradict the validity of (6) for uniformly elliptic, mixing environments. (See [4], page 21, for the definition of total ergodicity.)

Our attempts to address the validity of (6) in this last setting led to the tree tail estimates discussed in Theorem 2. Employing these bounds, we construct a counterexample to (6), in $d \geq 3$, with a stationary, uniformly elliptic and polynomially mixing environment.

THEOREM 3. *For $d \geq 3$, there is a probability space (with probability measure $\mathbb{P}$) that supports a stationary, uniformly elliptic and polynomially* *mixing family* $\omega = (\omega(x))_{x \in \mathbb{Z}^d}$, *such that for some constant* $c > 0$ *and* $\mathbb{P}$-*a.a. realizations of* $\omega$,

$$(7) \quad P_\omega^0 \left[ \liminf_{n \to \infty} \frac{X_n \cdot \vec{1}}{n} > c \right] > 0 \quad and \quad P_\omega^0 \left[ \liminf_{n \to \infty} \frac{X_n \cdot (-\vec{1})}{n} > c \right] > 0.$$

*Here,* $\vec{1} = e_1 + \cdots + e_d$.

We outline how we use the spanning forest constructed in Theorem 2 to obtain Theorem 3. The counterexample in [9] was based on constructing two disjoint directed trees in $\mathbb{Z}^2$ with opposite directions $z = \vec{1}$ and $z = -\vec{1}$, and adjusting the transition probabilities of the RWRE on each edge that belongs to one of the trees, so that the drift at $x$ toward the ancestor $a(x)$ increases as a function of $h(x)$. By appropriately choosing the rate at which the drift increases, one can ensure that the RWRE, when started on one of the trees, remains on it forever with positive probability, while progressing up its ancestral line. Because of this, the uniform ellipticity of the environment cannot be maintained.

When trying to restore uniform ellipticity to the environment, a natural idea is to add "insulation" around each of the directed trees. The insulation should allow one to specify a uniformly elliptic environment that, with positive probability, forever traps the walker near the tree. Of course, this implies that the insulation must grow as one progresses up the ancestral line. To leave room for two directed trees pointing in opposite directions to have nonoverlapping insulation, one needs for the trees not to be "too large." When quantifying the notion of "large" needed, one is naturally led to study the random variable $h(x)$ in (1).

The rest of the paper is organized as follows. Theorems 1 and 2 are proved in Section 2. In Section 3 we prune the forest obtained in Theorem 2 to make room for an independent copy of it with the direction $z$ reversed and then add insulation to be able to obtain uniform ellipticity of the environment of the RWRE later on. Geometric properties of insulated rays are investigated in Section 4. In Section 5 we equip each such insulated ray with an environment $\omega$ that traps the RWRE with positive probability. These environments are patched together in Section 6 to complete the proof of Theorem 3. After a short discussion of open problems in Section 7, we prove in the Appendix the mixing properties stated in Theorems 2 and 3.

We conclude the Introduction with some conventions and notation. The $p$-norm, $p \in [1, \infty]$ (on either $\mathbb{R}^d$ or $\mathbb{Z}^d$), will be denoted by $|\cdot|_p$. Most of the time, we will use the 1-norm, in which case we will drop the index 1 from $|\cdot|_1$. The metric $d(\cdot, \cdot)$ will always refer to $|\cdot|$, and $B(x, r)$ [resp. $B_\infty(x, r)$] denotes the closed $|\cdot|$ ball (resp. $|\cdot|_\infty$ ball) of center $x$ and radius $r$ in $\mathbb{Z}^d$. The collection of strictly positive integers will be denoted by $\mathbb{N}$ and the

cardinality of a set $A$ will be denoted by $\#A$. Throughout the paper, $c_i$, $i = 1, 2, 3, \ldots$, will denote strictly positive and finite constants that depend only on $d$ and $\beta$, where $\beta$ is introduced in (24).

**2. Spanning $\mathbb{Z}^d$ with short trees.** In this section we provide the proofs of Theorems 1 and 2. We begin with the easy proof of Theorem 1.

PROOF OF THEOREM 1.  Choose $c_1 > 0$ such that for all $n \geq 1$,
$$c_1 \#\{x \in \mathbb{Z}^d | |x| = n\} \leq n^{d-1}.$$

Since $(a(x))_{x \in \mathbb{Z}^d}$ is stationary, so is $(h(x))_{x \in \mathbb{Z}^d}$ and, therefore,

(8) $$n^{d-1} \mathbb{P}[h(0) \geq n] \geq c_1 \sum_{|x|=n} \mathbb{P}[h(x) \geq n] \qquad \text{for } n \geq 1.$$

If $x$ is an ancestor of $0$, then $h(x) \geq |x|$. Consequently, the right-hand side of (8) is at least
$$c_1 \sum_{|x|=n} \mathbb{P}[x \in \text{Ray}(0)] = c_1 \mathbb{E}[\#\{x \in \text{Ray}(0) | |x| = n\}] \geq c_1,$$

where the inequality holds since $\text{Ray}(0)$ contains at least one $x$ with $|x| = n$. The bound (2) follows.  □

The remainder of this section is devoted to demonstrating Theorem 2. We begin with the construction of the ancestral function $a$ referred to there. Fix $\gamma_d > 0$ such that for all $n \geq 1$,

(9) $$\#\{x \in \mathbb{N}^d | |x| = n\} \geq \gamma_d n^{d-1}.$$

Also, let $n_0 \in \mathbb{N}$ and $\theta_d$ be finite constants such that

(10) $$n_0^d \geq \theta_d \geq \frac{d^d}{\gamma_d}.$$

Let $L(x) > 1$, $x \in \mathbb{Z}^d$, be i.i.d. random variables whose distribution is atomless and satisfies

(11) $$P[L(0) > t] = \theta_d t^{-d} \qquad \text{for all } t \geq n_0.$$

We define, for each $t \geq 1$, the *umbrella*
$$U_t = \bigcup_{i=1}^{d} U_{i,t},$$

where
$$U_{i,t} = \{x \in [0,t]^d | x_i = 0, \text{ and } x_j > 0 \text{ for } j \neq i\}, \qquad i = 1, \ldots, d,$$

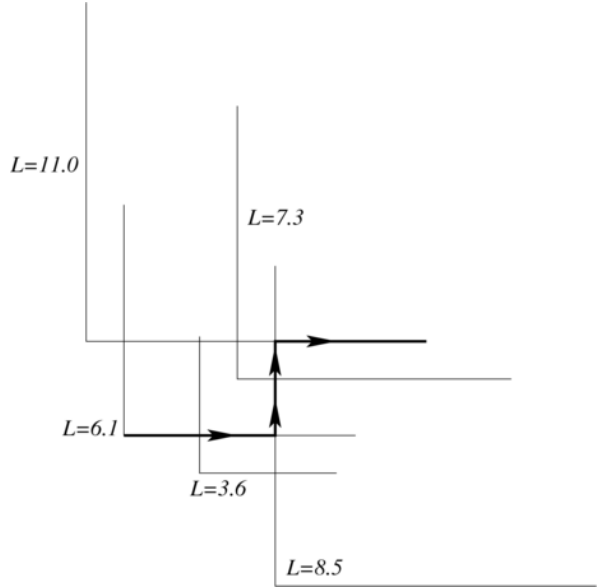

FIG. 2. *Water will follow the thick line if there are no umbrellas around larger than the ones shown.*

are the sides of the umbrella. Note that the umbrella $U_t$ contains exactly those points in $\mathbb{Z}^d$ through which one can enter the box $[1,t]^d \cap \mathbb{Z}^d$ by moving one step in one of the directions $e_1, \ldots, e_d$.

We next provide an informal description of the construction of the ancestral function $a$ that captures our way of thinking. [The reader can safely defer this discussion until after the precise definitions given in (12)–(14).] Imagine the ancestral line as being given by rainwater that always follows a directed nearest-neighbor flow on $\mathbb{Z}^d$. Rain is leaving each lattice point in all of the positive coordinate directions and is being deflected by the umbrellas $y + U_{L(y)}, y \in \mathbb{Z}^d$. The umbrellas will protect most of the points in the cubes $y + [1, L(y)]^d$ from the rain, as follows. Water that has reached a vertex $x$ is blocked from flowing to $x + e_i$ by any umbrella whose side $y + U_{i,L(y)}$ contains $x$. However, for every direction $e_1, \ldots, e_d$, there is an umbrella whose side $y + U_{i,L(y)} \ni x$ blocks that direction. This "battle" is lost by the direction $e_i$ for which the largest umbrella in $x$ blocking that direction is smallest among all directions, and the water will flow in this direction. (See Figure 2 for an illustration in $d = 2$. Note that the topological duality present for $d = 2$ is absent for $d \geq 3$.)

More precisely, we define, for all $i = 1, \ldots, d$ and all $x \in \mathbb{Z}^d$,

$$\lambda_i(x) = \sup_{y \in \mathbb{Z}^d \,:\, x \in y + U_{i,L(y)}} L(y), \tag{12}$$

which is the length of the largest umbrella whose $i$ side passes through $x$. Since $L(0) > 1$ a.s., we have $x \in x - e_j + U_{i,L(x-e_j)}$ a.s., for each $j \neq i$. Consequently, the set on the right-hand side of (12) is nonempty and $\lambda_i(x)$ is greater than 1. The following lemma implies that $\lambda_i(x)$ is also a.s. finite.

LEMMA 4. *There is a constant $c_3$, such that for all $i \in \{1, \ldots, d\}$ and all $t > n_0$,*

(13) $$\mathbb{P}[\lambda_i(0) > t] \leq c_3 t^{-1}.$$

PROOF. It suffices to show that (13) holds for $t$ large. Let $t_1 > n_0$ be chosen large enough so that $(1 - \theta_d t^{-d})^{t^d/\theta_d} > 1/3$ for $t > t_1$. Set

$$D_i^n = \{y \in \mathbb{Z}^d | |y|_\infty = n, y_i = 0, y_j < 0 \text{ for all } j \neq i\}.$$

By the definition of $\lambda_i(x)$ and the independence of $L(y), y \in \mathbb{Z}^d$, one obtains for $t > t_1$, that

$$\mathbb{P}[\lambda_i(0) \leq t] = \mathbb{P}[L(y) \leq t \text{ for all } y \in -U_{i,t}] \prod_{n>t} \mathbb{P}[L(y) < n \text{ for all } y \in D_i^n]$$

$$\overset{(11)}{\geq} (1 - \theta_d t^{-d})^{\lfloor t \rfloor^{d-1}} \prod_{n>t} (1 - \theta_d n^{-d})^{c_4 n^{d-2}}$$

$$\geq 3^{-\theta_d/t - \sum_{n>t} c_5 n^{-2}} \geq 3^{-c_6 t^{-1}} \geq 1 - c_3 t^{-1}$$

for appropriate constants $c_3, \ldots, c_6$. □

Since the distributions of $\lambda_i(x)$ and $L(x)$ are atomless, there is an a.s. unique $I(x) \in \{1, \ldots, d\}$, for which

$$\lambda_{I(x)}(x) = \min\{\lambda_i(x) | i = 1, \ldots, d\}.$$

This defines a direction with the smallest "protecting" umbrellas. Any umbrella passing through $x$ that is perpendicular to that direction will be penetrated at that site, in the sense that water will flow from $x$ in that direction. We now set, for $x \in \mathbb{Z}^d$,

(14) $$a(x) = x + e_{I(x)}.$$

Note that since $a$ is directed with $z = \vec{1}$, $a$ is an ancestral function; this is the ancestral function we will use to demonstrate Theorem 2. The edges $\{x, a(x)\}$, $x \in \mathbb{Z}^d$, that are "wetted by the rain" define a random forest of infinite trees that spans $\mathbb{Z}^d$ [as in Figure 1(b) for $d = 2$].

We still need to demonstrate the tail estimates and mixing properties in the statement of Theorem 2. As a first step, the next lemma bounds the probability that an umbrella, with side length at least $t$, has been penetrated at any given site.

LEMMA 5. *For some constant $c_7$ and all $i \in \{1, \ldots, d\}$, $t > n_0$ and $z \in U_{i,t}$,*

(15) $$\mathbb{P}[I(z) = i, L(0) > t] \leq c_7 t^{-2d+1}.$$

PROOF. Denote by $\mathcal{A}_i(z)$ the $\sigma$-field generated by the random variables $L(u)$, with $u_i = z_i$ and $u_j \neq z_j$ for all $j \neq i$. Note that $\mathcal{A}_i(z)$, $i = 1, \ldots, d$, are independent. Moreover, $\lambda_i(z)$ is measurable with respect to $\mathcal{A}_i(z)$, and $L(0)$ is measurable with respect to $\mathcal{A}_i(z)$ if $z \in U_{i,n}$. Therefore, for $i = 1, \ldots, d$ and $z \in U_{i,n}$,

$$\mathbb{P}[I(z) = i, L(0) > t] = \mathbb{E}\left[\mathbb{P}\left[\lambda_i(z) < \min_{j \neq i} \lambda_j(z) \Big| \mathcal{A}_i(z)\right]; L(0) > t\right]$$
$$= \mathbb{E}\left[\prod_{j \neq i} \mathbb{P}[\lambda_i(z) < \lambda_j(z) | \mathcal{A}_i(z)]; L(0) > t\right].$$

By (13) (for the first inequality), and by (11) and the independence of $(L(x))_x$ (for the second inequality), this is

$$\leq \mathbb{E}[(c_3 \lambda_i(z)^{-1})^{d-1}; L(0) > t] \leq c_7 t^{-(d-1)-d},$$

because $L(0) > t$ implies $\lambda_i(z) > t$. □

We can now demonstrate Theorem 2.

PROOF OF THEOREM 2. It remains to show that the ancestral function $(a(x))_x$ constructed above is polynomially mixing of order 1 and satisfies the bound in (3). Here, we demonstrate (3); the demonstration of polynomial mixing is deferred to Lemma A.1 in the Appendix. (The lemma in fact deals with a slightly stronger notion of mixing.)

We begin the proof of (3) by introducing some notation. Denote by

$$S_m^n = \{x \in \mathbb{Z}^d | m \leq x \cdot \vec{1} \leq n\}, \qquad m, n \in \mathbb{Z},$$

the slab bounded by the hyperplanes perpendicular to $\vec{1}$ and passing through $m\vec{1}$ and $n\vec{1}$. (This slab is empty if $m > n$.) Also, set

$$S_m^{n,+} = S_m^n \cap \mathbb{N}^d, \qquad S_m^{n,-} = S_m^n \cap -\mathbb{N}^d.$$

We define the random variables

(16) $\quad M_n = \sup\{m \in \{n_0, \ldots, n\} | \exists x \in S_{-m}^{-m,-}$ with $L(x) > m\}, \qquad n \geq 1,$

where we set $M_n = n_0 - 1$ if the set in (16) is empty.

We will show that

(17) $\qquad \mathbb{P}[h(0) > m, M_n = m] \leq c_8 n^{-d}, \qquad m = n_0 - 1, \ldots, n.$

This implies

$$\mathbb{P}[h(0) > n] \leq \sum_{m=n_0-1}^{n} P[h(0) > m, M_n = m] \leq c_9 n^{-d+1},$$

from which (3) follows. The proof of (17) is divided into the degenerate case, $m = n_0 - 1$, and the general case, $m = n_0, \ldots, n$, which involves considerably more work.

We first consider the case $m = n_0 - 1$. Then, there is no umbrella $y + U_{L(y)}$, with $n_0 \leq |y| \leq n$, that protects the origin. Therefore,

$$\mathbb{P}[h(0) > m, M_n = n_0 - 1]$$
$$\leq \mathbb{P}[M_n = n_0 - 1]$$
$$= \mathbb{P}[\text{for all } m = n_0, \ldots, n \text{ and } y \in S_{-m}^{-m,-}, L(y) \leq m]$$
$$\stackrel{(11)}{=} \prod_{m=n_0}^{n} \left(1 - \frac{\theta_d}{m^d}\right)^{\#S_{-m}^{-m,-}} \stackrel{(9)}{\leq} \prod_{m=n_0}^{n} \left(1 - \frac{\theta_d}{m^d}\right)^{(m^d/\theta_d)(\gamma_d \theta_d/m)}.$$

Since $(1 - x^{-1})^x \leq e^{-1}$ for $x \geq 1$, the last expression is, by (10) and (11), at most

$$(18) \qquad \exp\left(-d^d \sum_{m=n_0}^{n} m^{-1}\right) \leq \exp\left(-d^d \int_{n_0}^{n} x^{-1} \, dx\right) \leq c_{10} n^{-d}$$

for appropriate $c_{10}$. This implies (17) for $m = n_0 - 1$ if one takes $c_8 \geq c_{10}$.

We now demonstrate (17) for the general case $m = n_0, \ldots, n$. We will employ the events

$$A_m^n(x, r) = \{L(y) \leq -y \cdot \vec{1} + r \text{ for all } y \in S_m^n \cap (x - \mathbb{N}^d)\},$$
$$m, n, r \in \mathbb{Z}, x \in \mathbb{Z}^d.$$

Since the proof is long, we break it into three parts. The first part consists of showing

$$(19) \quad \begin{aligned} &\mathbb{P}[h(0) > m, M_n = m] \\ &\leq \sum_{s=1}^{d} \mathbb{E}\left[\#\{x \in S_m^{m,+} | h(x) > m\}; L(0) > m; A_{m-n}^{-1}\left(\left\lfloor \frac{m}{d} \right\rfloor e_s, m\right)\right]. \end{aligned}$$

In the next part, we will decouple the events that appear in this expectation, whose probabilities we will then compute.

To prove (19), first note that by definition (16),

$$\mathbb{P}[h(0) > m, M_n = m]$$
$$= P[h(0) > m, L(x) > m \text{ for some } x \in S_{-m}^{-m,-},$$

$$L(y) \leq -y \cdot \vec{1} \text{ for all } y \in S_{-n}^{-m-1,-}]$$
$$\leq \sum_{x \in S_{-m}^{-m,-}} \mathbb{P}[h(0) > m, L(x) > m, A_{-n}^{-m-1}(0,0)].$$

By stationarity, this equals

$$\sum_{x \in S_{-m}^{-m,-}} \mathbb{P}[h(-x) > m, L(0) > m, A_{m-n}^{-1}(-x, m)]$$

$$= \sum_{x \in S_m^{m,+}} \mathbb{P}[h(x) > m, L(0) > m, A_{m-n}^{-1}(x, m)]$$

$$\leq \sum_{s=1}^{d} \sum_{\substack{x \in S_m^{m,+} \\ x_s = |x|_\infty}} \mathbb{P}[h(x) > m, L(0) > m, A_{m-n}^{-1}(x, m)].$$

Observe that $(m/d)e_s \leq x$ coordinatewise whenever $x \in S_m^{m,+}$ and $x_s = |x|_\infty$. For such $x$, $\lfloor m/d \rfloor e_s - \mathbb{N}^d \subseteq x - \mathbb{N}^d$. (See Figure 3 for an illustration.) Consequently, the above double sum is at most

$$\sum_{s=1}^{d} \sum_{\substack{x \in S_m^{m,+} \\ x_s = |x|_\infty}} \mathbb{P}\left[h(x) > m, L(0) > m, A_{m-n}^{-1}\left(\left\lfloor \frac{m}{d} \right\rfloor e_s, m\right)\right],$$

which yields (19).

We next perform the decoupling previously mentioned, which leads to (21) below. If $x \in S_m^{m,+}$ and $h(x) > m$, then $\text{Tree}(x)$ is not contained in the cube $[1, m]^d$ (see Figure 3). That is, $\text{Tree}(x)$ must possess at least one branch that penetrates the umbrella $U_m$ that "protects" the cube. Consequently,

$$\#\{x \in S_m^{m,+} | h(x) > m\} \leq \#\left(\bigcup_{z : z \in U_{I(z),m}} \text{Ray}(z) \cap S_m^{m,+}\right)$$

$$\leq \sum_{z : z \in U_{I(z),m}} \#(\text{Ray}(z) \cap S_m^m) = \sum_{i=1}^{d} \sum_{z \in U_{i,m}} \mathbf{1}_{I(z)=i}.$$

In the last step, we used the fact that rays are directed (in the direction $\vec{1}$) and therefore can intersect the hyperplane $S_m^m$ at exactly one site. Substituting this into (19) yields

(20)
$$\mathbb{P}[h(0) > m, M_n = m]$$
$$\leq \sum_{s=1}^{d} \sum_{i=1}^{d} \sum_{z \in U_{i,m}} \mathbb{P}\left[I(z) = i, L(0) > m, A_{m-n}^{-1}\left(\left\lfloor \frac{m}{d} \right\rfloor e_s, m\right)\right].$$

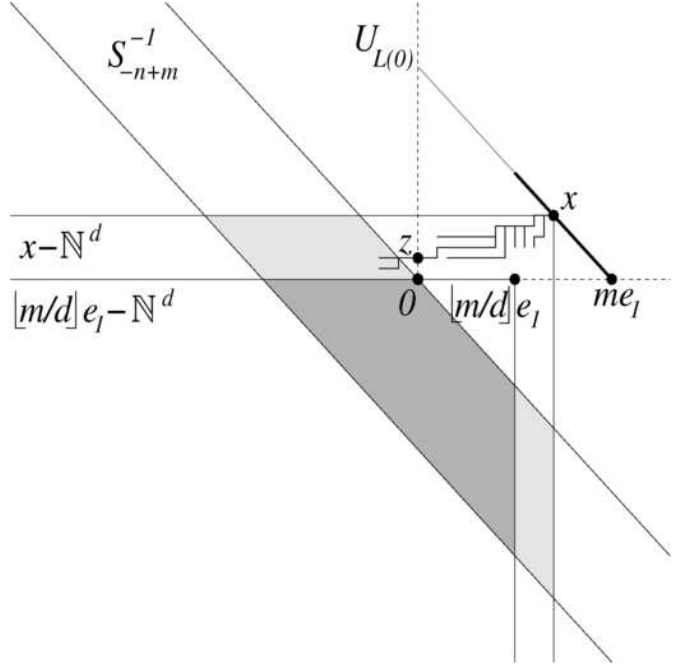

Fig. 3. *The dark region* $(\lfloor m/d \rfloor e_1 - \mathbb{N}^d) \cap S_{m-n}^{-1}$ *is always included in the lightly shaded region* $(x - \mathbb{N}^d) \cap S_{m-n}^{-1}$, *regardless of where $x$ is located on the bold line. The condition $h(x) > m$ implies that* Tree$(x)$ *extends past the box on whose diagonal $x$ is located.*

To get (17), we would like the first two events inside the last probability to be independent of $A_{m-n}^{-1}(\lfloor m/d \rfloor e_s, m)$. Unfortunately, this is not quite true, but will be true if we replace it by the larger event $B_{m-n}^{-1}(\lfloor m/d \rfloor e_s, m)$, where

$$B_{m-n}^{-1}(\lfloor m/d \rfloor e_s, m)$$
$$= \Bigg\{ L(y) \leq m - y \cdot \vec{1}$$
$$\text{for all } y \in \left( S_{m-n}^{-1} \cap \left( \left\lfloor \frac{m}{d} \right\rfloor e_s - \mathbb{N}^d \right) \right) \Big\backslash \bigcup_{j=1}^{d} \{x | x_j = z_j\} \Bigg\}.$$

Indeed, for $z \in U_{i,m}$, $I(z)$ and $L(0)$ are measurable with respect to the $\sigma$-field generated by the random variables $L(x)$, where $x$ has at least one coordinate in common with $z$, whereas $B_{m-n}^{-1}(\lfloor m/d \rfloor e_s, m)$ is independent of $L(x)$ for such $x$. Therefore, we obtain

$$\mathbb{P}[h(0) > m, M_n = m] \tag{21}$$

$$\leq \sum_{s,i=1}^{d} \sum_{z \in U_{i,m}} \mathbb{P}[I(z) = i, L(0) > m]\mathbb{P}[B_{m-n}^{-1}(\lfloor m/d \rfloor e_s, m)].$$

As the last step in showing (17), we bound the probabilities that appear on the right-hand side of (21). By Lemma 5, (11) and independence of $(L(x))_x$, the right-hand side of (21) is at most

(22)
$$\sum_{s,i=1}^{d} \sum_{z \in U_{i,m}} c_7 m^{-2d+1}$$
$$\times \prod_{k=2}^{n-m} \left(1 - \frac{\theta_d}{(m+k)^d}\right)^{\#((S_{-k}^{-k} \cap (\lfloor m/d \rfloor e_s - \mathbb{N}^d)) \setminus \bigcup_{j=1}^{d} \{x | x_j = z_j\})},$$

where we have dropped the factor for $k = 1$.

We proceed to simplify (22) by estimating the cardinality of the sets in the exponent of the above product. For all $m \geq 1, k \geq 2$ and $s = 1, \ldots, d$,

$$\#\left(S_{-k}^{-k} \cap \left(\left\lfloor \frac{m}{d} \right\rfloor e_s - \mathbb{N}^d\right)\right) \overset{(9)}{\geq} \gamma_d \left(\left\lfloor \frac{m}{d} \right\rfloor + k\right)^{d-1} \geq \gamma_d \left(\frac{m + (k-1)d}{d}\right)^{d-1}$$
$$\geq \gamma_d d^{-d+1}(m+k)^{d-1} \overset{(10)}{\geq} \frac{d}{\theta_d}(m+k)^{d-1}.$$

[For the first inequality, note that the set on the left-hand side consists of those $z \in \mathbb{Z}^d$ with $z < \lfloor m/d \rfloor e_s$ coordinatewise and $d(z, \lfloor m/d \rfloor e_s) = \lfloor m/d \rfloor + k$.] Similarly,

$$\#\left(S_{-k}^{-k} \cap \left(\left\lfloor \frac{m}{d} \right\rfloor e_s - \mathbb{N}^d\right) \cap \bigcup_{j=1}^{d} \{x | x_j = z_j\}\right) \leq c_{11} \left(\frac{m}{d} + k\right)^{d-2}.$$

The expression in (22) is consequently at most

(23)
$$\sum_{s,i=1}^{d} \sum_{z \in U_{i,m}} c_7 m^{-2d+1} \prod_{k=2}^{n-m} \left(1 - \frac{\theta_d}{(m+k)^d}\right)^{d(m+k)^{d-1}/\theta_d}$$
$$\times \prod_{k \geq 2} \left(1 - \frac{\theta_d}{(m+k)^d}\right)^{-c_{11}(m/d+k)^{d-2}}.$$

One easily checks that the infinite product in (23) is less than a constant independent of $m$ because of the difference of 2 in the exponents $d$ and $d-2$ of $k$. [Recall that $\theta_d < (m+k)^d$ for $k \geq 2$, by (10).] Applying $(1-x^{-1})^x \leq e^{-1}$ for $x \geq 1$ to the terms of the finite product, the right-hand side of (23) is at

most

$$d^2(\#U_{1,m})c_{12}m^{-2d+1}\exp\left(-d\sum_{k=2}^{n-m}(m+k)^{-1}\right)$$
$$\leq c_{13}m^{d-1}m^{-2d+1}\exp\left(-d\sum_{k=m+2}^{n}k^{-1}\right)$$
$$\leq c_{13}m^{-d}\exp\left(-d\int_{m+2}^{n}t^{-1}\,dt\right) = c_{13}\left(\frac{m+2}{m}\right)^d n^{-d} \leq c_{14}n^{-d}.$$

This bounds $P[h(0) > m, M_n = m]$ for $m = n_0, \ldots, n$. Together with (18), this completes the proof of (17) if one chooses $c_8 = \max(c_{14}, c_{10})$. □

**3. Pruning and insulating trees.** For the remainder of the paper, we will consider arbitrary ancestral functions $(a(x))_{x \in \mathbb{Z}^d}$ that satisfy the statement of Theorem 2. Eventually, when studying mixing properties (in the proof of Theorem 3 and in the Appendix), we will need to use the explicit construction of ancestral functions provided in Section 2.

In this section we *prune* the trees constructed in Theorem 2, with an eye toward the construction of an environment for the RWRE. This will allow us to build wide enough *channels* to trap the random walker and direct it in certain chosen directions. Throughout the remainder of the paper, we will assume $d \geq 3$ and employ a constant $\beta$ that satisfies

(24) $$0 < \beta < \frac{d-2}{2d}.$$

(This bound is needed for the construction of the environment. To guarantee its mixing properties, we will eventually take $\beta$ small enough so that the conclusion of Lemma A.5 in the Appendix holds.)

For each $y \in \mathbb{Z}^d$, we consider the ball $B(y, h(y)^\beta)$, where $h$ is defined in (1). (Such balls will serve as the "insulation" alluded to in the Introduction.) Any given $x \in \mathbb{Z}^d$ may be covered by a number of balls: we set

(25) $$H(x) = \sup\{h(y) | x \in B(y, h(y)^\beta)\}.$$

We first estimate the tail behavior of $H(0)$.

LEMMA 6. *For appropriate $c_{17}$,*

(26) $$\limsup_{n \to \infty} n^{(1-\beta)d-1} \mathbb{P}[H(0) \geq n] \leq c_{17} < \infty.$$

PROOF. For $t = n^\beta$ chosen large enough,

$$\mathbb{P}[H(0)^\beta \geq t] = \mathbb{P}[0 \in B(y, h(y)^\beta) \text{ and } h(y)^\beta \geq t, \text{ for some } y \in \mathbb{Z}^d]$$

$$\leq \sum_{m \geq 0} \sum_{|y|=m} \mathbb{P}[h(y)^\beta \geq m \vee t]$$

$$\leq \sum_{|y| \leq \lfloor t \rfloor} \mathbb{P}[h(0) \geq t^{1/\beta}] + \sum_{m > \lfloor t \rfloor} \sum_{|y|=m} \mathbb{P}[h(0) \geq m^{1/\beta}].$$

By (3) of Theorem 2, this is

$$\leq c_{15} t^d t^{(1-d)/\beta} + \sum_{m > \lfloor t \rfloor} c_{16} m^{d-1} m^{(1-d)/\beta} \leq c_{17} t^{d+(1-d)/\beta}$$

for appropriate $c_{15}, c_{16}, c_{17}$, since $d - 1 > \beta d$ by (24). This implies (26). □

Recall that Theorem 2 guarantees the existence of a directed stationary forest with certain specific properties. By relabeling the coordinate axes, we can choose the directions in which the rays of this forest grow, that is, we can specify $z \in \{\pm 1\}^d$ such that, for each $x \in \mathbb{Z}^d$, $a(x) - x \in \{z_i e_i | i = 1, \ldots, d\}$. We may therefore assume that there exist, on the same probability space, two independent directed forests, $F_i = \{\{x, a_i(x)\} | x \in \mathbb{Z}^d\}$, $i = 1, 2$, with ancestral functions $a_i$, such that $a_1(x) - x = e_j$ for some $j = 1, \ldots, d$ and $a_2(x) - x = -e_j$ for some $j = 1, \ldots, d$. These two forests "grow" in opposite directions $\vec{1}$ and $-\vec{1}$. We define the corresponding functions $h_i$ and $H_i$ in the same way as above, using their respective forests $F_i$.

We proceed to "prune" the trees in both forests in such a way that each pruned forest will consist solely of infinite trees, and these forests will be "well separated." (The forests will no longer span $\mathbb{Z}^d$.) To this end, we define, for $(i, j) = (1, 2), (2, 1)$,

$$\widetilde{\mathbb{T}}_i = \{x \in \mathbb{Z}^d | h_i(x) > H_j(y) \text{ for } y \in B(x, h_i(x)^\beta)\}.$$

We prune the original forests $F_i$ by removing the vertices $\mathbb{Z}^d \setminus \widetilde{\mathbb{T}}_i$. This will split any tree in such a forest into a number of finite and at most one infinite piece. After removing the finite branches, we are left with the set of (a priori, possibly empty) directed infinite pruned trees

$$\mathbb{T}_i = \{x \in \widetilde{\mathbb{T}}_i | a_i^n(x) \in \widetilde{\mathbb{T}}_i \text{ for all } n \geq 0\}.$$

So that the transition probabilities we are going to construct are uniformly elliptic, we *insulate* the forests $\mathbb{T}_i$ by defining, for $i = 1, 2$,

(27) $$B_i = \bigcup_{x \in \mathbb{T}_i} B(x, h_i(x)^\beta).$$

The next proposition shows that the sets $B_1$ and $B_2$ are disjoint and not empty.

PROPOSITION 7. *The sets $B_1$ and $B_2$ are a.s. disjoint. There exist a.s. $N_i \in \mathbb{N}$, $i = 1, 2$, such that $a_i^n(0) \in \mathbb{T}_i$ for all $n \geq N_i$. In particular, $B_i \neq \varnothing$ a.s.*

PROOF. To prove the disjointness of $B_1$ and $B_2$, assume instead that $x \in B_1 \cap B_2$. Then there exist $x_1 \in \mathbb{T}_1$ and $x_2 \in \mathbb{T}_2$ so that $x \in B(x_i, h_i(x_i)^\beta)$ for $i = 1, 2$. Since $x \in B(x_1, h_1(x_1)^\beta)$, one has $H_1(x) \geq h_1(x_1)$. Moreover, since $x_2 \in \widetilde{\mathbb{T}}_2$ and $x \in B(x_2, h_2(x_2)^\beta)$, one has $H_1(x) < h_2(x_2)$ by the definition of $\mathbb{T}_2$. Consequently, $h_2(x_2) > h_1(x_1)$. Analogously, one obtains $h_1(x_1) > h_2(x_2)$, which is a contradiction and proves $B_1 \cap B_2 = \varnothing$.

We will next show that

$$(28) \quad \limsup_{k \to \infty} k^{(1-2\beta)d-2} \mathbb{P}[a_i^n(0) \notin \widetilde{\mathbb{T}}_i \text{ for some } n \geq k] < \infty, \qquad i = 1, 2,$$

which implies the second claim since, by (24), $(1 - 2\beta)d - 2 > 0$. To demonstrate (28), let $k \geq 0$ and $(i, j) = (1, 2)$ or $(i, j) = (2, 1)$. Then

$$\mathbb{P}[a_i^n(0) \notin \widetilde{\mathbb{T}}_i \text{ for some } n \geq k]$$
$$(29) \quad \leq \sum_{n \geq k} \mathbb{P}[H_j(y) \geq h_i(a_i^n(0)) \text{ for some } y \in B(a_i^n(0), h_i(a_i^n(0))^\beta)]$$
$$\leq \sum_{n \geq k} \mathbb{E}\left[\sum_{y \in B(a_i^n(0), h_i(a_i^n(0))^\beta)} \mathbb{P}[H_j(y) \geq h_i(a_i^n(0)) | \sigma(a_i(x), x \in \mathbb{Z}^d)]\right].$$

The number of terms in the inner sum is bounded above by $c_{18}(h_i(a_i^n(0))^{\beta d}$ for appropriate $c_{18}$. Moreover, $H_j$ is measurable with respect to $\sigma(a_j(x), x \in \mathbb{Z}^d)$, which is independent of $\sigma(a_i(x), x \in \mathbb{Z}^d)$. Hence, we can use Lemma 6 to estimate each such term from above and obtain that (29) is at most

$$\sum_{n \geq k} \mathbb{E}[c_{19} h_i(a_i^n(0))^{1-(1-2\beta)d}] \leq \sum_{n \geq k} c_{19} n^{1-(1-2\beta)d}$$

for appropriate $c_{19}$, since $h_i(a_i^n(0)) \geq n$. Since $1 - (1 - 2\beta)d < -1$ by (24), inequality (28) follows. □

Since they are subsets of $F_i$, the sets $\Pi_i = \{(x, a_i(x)) | x \in \mathbb{T}_i\}$, $i = 1, 2$, are also forests. By Proposition 7, they almost surely contain infinite rays that point in opposite directions and they do not have any vertices in common. Moreover, the set of immediate ancestors of vertices in $\mathbb{T}_i$ is contained in $B_i$, because, for any $x \in \mathbb{T}_i$, $h_i(x) \geq 1$ [since $H_j(x) \geq 0$] and hence $a_i(x) \in B(x, h_i(x)^\beta)$. Consequently, no vertex in $\Pi_1$ is connected to a vertex in $\Pi_2$ and so $\Pi_1 \cup \Pi_2$ is also a forest (although it does not span all of $\mathbb{Z}^d$). With a slight abuse of notation, we say that the ancestral function of $\Pi_1 \cup \Pi_2$ is given by

$$\alpha(x) = \begin{cases} a_1(x), & \text{if } x \in \mathbb{T}_1, \\ a_2(x), & \text{if } x \in \mathbb{T}_2, \end{cases} \tag{30}$$

where $\alpha(x)$ is defined only for $x \in \mathbb{T}_1 \cup \mathbb{T}_2$.

**4. Geometry of insulated rays.** In this section we introduce terminology and provide estimates that we will need when we analyze the RWRE environment in Sections 5 and 6. For $i = 1, 2$, let $\partial \mathbb{T}_i = \{z \in \mathbb{T}_i | z \neq \alpha(x) \text{ for all } x \in \mathbb{T}_i\}$ denote the set of *leaves* of the infinite pruned tree $\mathbb{T}_i$. By (27),

$$B_i = \bigcup_{z \in \partial \mathbb{T}_i} \bigcup_{n \geq 0} B(\alpha^n(z), (h_i(\alpha^n(z)))^\beta), \qquad i = 1, 2.$$

Instead of $B_i$, we will work with the somewhat simpler sets (see Figure 4)

$$C_i = \bigcup_{z \in \partial \mathbb{T}_i} \operatorname{InsRay}(z), \tag{31}$$

where

$$\operatorname{InsRay}(z) = \bigcup_{n \geq 0} B(\alpha^n(z), n^\beta) \tag{32}$$

is the *insulated ray emanating from* $z \in \partial \mathbb{T}_i$. [Since $\partial \mathbb{T}_1$ and $\partial \mathbb{T}_2$ are disjoint, there is no need to index $\operatorname{InsRay}(z)$ or $\operatorname{Ray}(z)$ with $i$.] Because $h_i(\alpha^n(z)) \geq n$, $C_i \subseteq B_i$. In particular, since $B_1$ and $B_2$ are disjoint by Proposition 7, so are $C_1$ and $C_2$.

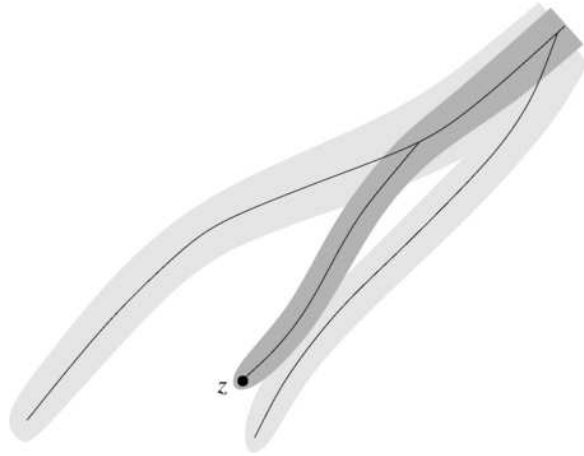

FIG. 4. $C_i$ is shaded and $\operatorname{InsRay}(z)$ is darkly shaded.

For $z \in \partial \mathbb{T}_1 \cup \partial \mathbb{T}_2$ and $x \in \mathbb{Z}^d$, we define two quantities $u_z(x)$ and $v_z(x)$, each measuring the "effective" insulation at $x$ in a slightly different way. We set

$$(33) \quad u_z(x) = d(x, \mathrm{InsRay}(z)^c) \quad \text{and} \quad v_z(x) = \sup_{n \geq 0}(n^\beta - |x - \alpha^n(z)|).$$

One can check that

$$(34) \qquad v_z(x) \leq \sup_{n \geq 0} d(x, B(\alpha^n(z), n^\beta)^c) \leq u_z(x).$$

Also, let $n_z(x)$ be the largest value of $n$ at which the supremum in (33) is attained, that is,

$$(35) \qquad n_z(x) = \max\{n \geq 0 : n^\beta - |x - \alpha^n(z)| = v_z(x)\}.$$

Since $u_z(x) < \infty$ and $\beta < 1$, $n_z(x) < \infty$, a.s. Also,

$$(36) \qquad |v_z(x) - v_z(x+e)| \leq 1$$

for all $x, e \in \mathbb{Z}^d$ with $|e| = 1$; this follows from

$$\begin{aligned}
v_z(x) &= n_z(x)^\beta - |x - \alpha^{n_z(x)}(z)| \\
&\leq n_z(x)^\beta - |x + e - \alpha^{n_z(x)}(z)| + |e| \\
&\leq v_z(x+e) + 1.
\end{aligned}$$

We will need the following estimates that involve $u_z(x)$ and $x - z$ in Sections 5 and 6.

LEMMA 8. *For appropriate $c_{20}$, all $z \in \partial \mathbb{T}_i$, $i = 1, 2$, and all $x \in \mathrm{InsRay}(z)$,*

$$(37) \qquad |x - z| \leq 2H_i(x)$$

*and*

$$(38) \qquad u_z(x) \leq c_{20}((-1)^{i+1}(x - z) \cdot \vec{1})^\beta.$$

*Consequently, for appropriate $c_{21}$,*

$$(39) \qquad u_z(x) \leq c_{21} H_i(x)^\beta.$$

PROOF. Inequality (39) follows from (37) and (38). For the proof of (37), recall that since $x \in \mathrm{InsRay}(z)$ for $z \in \partial \mathbb{T}_i$, $|x - a_i^m(z)| \leq m^\beta$ for some $m \geq 0$; in particular, $H_i(x) \geq h_i(a_i^m(z)) \geq m$. Therefore,

$$|x - z| \leq |x - a_i^m(z)| + |a_i^m(z) - z| \leq m^\beta + m \leq 2m \leq 2H_i(x).$$

The argument for (38) is longer. Set $n = (-1)^{i+1}(x - z) \cdot \vec{1}$; $n \geq 0$ since $\mathrm{Ray}(z)$ is directed. Moreover, we may assume $n \geq 1$, because $n = 0$ implies that $x = z$, in which case $u_z(x) = 0$ and (38) is trivial.

We introduce a vector $w \neq 0$, as follows. For $x \notin \text{Ray}(z)$, one has $x \neq \alpha^n(z)$, in which case we set $w = x - \alpha^n(z)$. Then
$$w \cdot \vec{1} = (x - z) \cdot \vec{1} + (z - \alpha^n(z)) \cdot \vec{1} = (-1)^{i+1}(n - n) = 0,$$
since $(z - \alpha^n(z)) \cdot \vec{1} = (-1)^i n$. For $x \in \text{Ray}(z)$, one has $x = \alpha^n(z)$; we then set $w = e_1 - e_2 \neq 0$, which is also orthogonal to $\vec{1}$.

We choose $c_{20}$ large enough so that

(40) $\qquad c_{22} = ((c_{20} - 2)d^{-2})^{1/\beta} \qquad \text{satisfies } c_{22}^{-\beta}(c_{22} - 1) > \sqrt{d}$

and set
$$y = x + \left\lfloor \frac{c_{20}n^\beta}{|w|} \right\rfloor w = \alpha^n(z) + \left(\mathbf{1}\{x \notin \text{Ray}(z)\} + \left\lfloor \frac{c_{20}n^\beta}{|w|} \right\rfloor\right)w.$$

To demonstrate (38), it is enough to show $y \notin \text{InsRay}(z)$, since then $u_z(x) \leq |x - y| \leq c_{20}n^\beta$.

We argue by contradiction and assume that $y \in \text{InsRay}(z)$. We will show that there exists an $m$ such that both

(41) $\qquad\qquad\qquad\qquad m \geq c_{22}n$

and

(42) $\qquad\qquad\qquad\qquad m \leq \left(\frac{c_{22}\sqrt{d}}{c_{22} - 1}\right)^{1/(1-\beta)}$

must hold. This is not possible because of our choice of $c_{22}$ in (40) and $n \geq 1$.

We choose $m \geq 0$ so that

(43)
$$m^\beta \geq |y - \alpha^m(z)|$$
$$= \left|\alpha^n(z) - \alpha^m(z) + \left(\mathbf{1}\{x \notin \text{Ray}(z)\} + \left\lfloor \frac{c_{20}n^\beta}{|w|} \right\rfloor\right)w\right|$$

and show that $m$ satisfies (41) and (42). Since $\text{Ray}(z)$ is directed, the coordinates of $\alpha^n(z) - \alpha^m(z)$ all have the same sign. On the other hand, $|w_j| \geq d^{-1}|w|$ must hold for at least one of the coordinates $w_j$ of $w$. There is also at least one other coordinate $w_k$ of $w$ with sign opposite to that of $w_j$ and with $|w_k| \geq d^{-2}|w|$, because $w \cdot \vec{1} = 0$. Therefore, either $w_j$ or $w_k$ has sign that is the opposite of that of the corresponding coordinate of $\alpha^n(z) - \alpha^m(z)$ and has absolute value at least $d^{-2}|w|$. So, by (43),

$$m^\beta \geq \left(\mathbf{1}\{x \notin \text{Ray}(z)\} + \left\lfloor \frac{c_{20}n^\beta}{|w|} \right\rfloor\right)d^{-2}|w|$$
$$\geq \left(\frac{c_{20}n^\beta}{|w|} - \mathbf{1}\{x \in \text{Ray}(z)\}\right)d^{-2}|w|$$
$$\geq (c_{20} - 2)d^{-2}n^\beta.$$

Because of our choice of $c_{22}$ in (40), this implies (41).

We still need to show (42), which we do by bounding $|y - \alpha^m(z)|$ from below. One has the string of inequalities

$$\begin{aligned}
|y - \alpha^m(z)| &\geq |y - \alpha^m(z)|_2 \\
&\geq |(y - \alpha^m(z)) \cdot \vec{1}|/\sqrt{d} \\
&= |(\alpha^n(z) - \alpha^m(z)) \cdot \vec{1}|/\sqrt{d} \\
&= |m - n|/\sqrt{d} \\
&\geq (1 - 1/c_{22})m/\sqrt{d},
\end{aligned}$$

where the second inequality follows from Cauchy–Schwarz and the last inequality follows from (41) and $c_{22} > 1$. Along with (43) this implies (42) and completes the proof of (38). □

**5. Environment attached to an insulated ray.** In this section we construct, for every insulated ray $\text{InsRay}(z)$, $z \in \partial \mathbb{T}_1 \cup \partial \mathbb{T}_2$, a uniformly elliptic environment that forever traps the walk $\{X_n\}$ inside $\text{InsRay}(z)$ with positive probability. This is achieved by both "pushing" the walk toward $\text{Ray}(z)$ and in a direction parallel to $\text{Ray}(z)$ in which $\text{InsRay}(z)$ widens. Proposition 7 is the main result of the section; most of the work is done in Lemmas 10 and 11.

These two directions are determined as follows, for any $x \in \text{InsRay}(z)$. The parallel motion of the walk consists of jumping from $x$ to $x + r_z(x)$, where

(44) $$r_z(x) = \alpha^{n_z(x)+1}(z) - \alpha^{n_z(x)}(z)$$

and $n_z(x)$ is given by (35). To define the second direction $s_z(x)$, note that since

(45) $$x \in \text{Ray}(z) \qquad \text{iff } x = \alpha^{n_z(x)}(z),$$

there exists, for any $x \notin \text{Ray}(z)$, a (deterministically chosen) unit vector $s_z(x) \in \mathbb{Z}^d$ so that

(46) $$|x + s_z(x) - \alpha^{n_z(x)}(z)| = |x - \alpha^{n_z(x)}(z)| - 1;$$

moving from $x$ to $x + s_z(x)$ takes the walk closer to $\alpha^{n_z(x)}(z)$. For $x \in \text{Ray}(z)$, the choice of the unit vector $s_z(x)$ is arbitrary. (See Figure 5 for an illustration of both motions, for $z \in B_1$.)

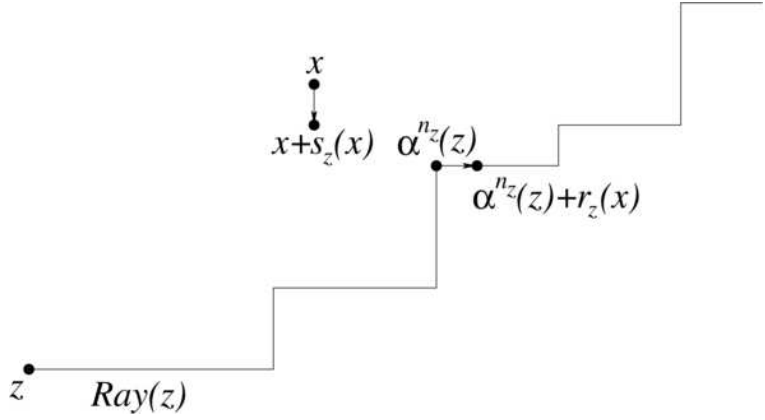

Fig. 5.

For $z \in \partial \mathbb{T}_1 \cup \partial \mathbb{T}_2$, we define the *environment* $\omega^z$ *attached to* $\mathrm{Ray}(z)$ by setting, for $e \in \mathbb{Z}^d$ with $|e| = 1$,

$$(47) \quad \omega^z(x, x+e) = \begin{cases} \dfrac{3}{4}\mathbf{1}_{e=r_z(x)} + \dfrac{1}{5}\mathbf{1}_{e=s_z(x)} + \dfrac{\mathbf{1}_{e \notin \{r_z(x), s_z(x)\}}}{20(2d - 1 - \mathbf{1}_{r_z(x) \neq s_z(x)})}, \\ \qquad \text{if } x \in \mathrm{InsRay}(z), \\ \dfrac{1}{2d}, \qquad \text{if } x \notin \mathrm{InsRay}(z). \end{cases}$$

That is, when $x \in \mathrm{InsRay}(z)$, the walk moves with probabilities $3/4$ and $1/5$ in the directions $r_z(x)$ and $s_z(x)$, respectively (if the directions are different), and uniformly chooses one of the other directions with the remaining probability; when $x \notin \mathrm{InsRay}(z)$, the motion of the walk is symmetric.

The environment $\omega^z \in \Omega$ is uniformly elliptic since for $x, e \in \mathbb{Z}^d$, with $|e| = 1$, and for $\kappa = (20(2d-1))^{-1}$,

$$(48) \qquad \omega^z(x, x+e) \geq \kappa.$$

As mentioned in the beginning of the section, $\omega^z$ has been constructed so that, for $X_0 \in \mathrm{InsRay}(z)$,

$$T_z = \inf\{n \geq 0 : X_n \notin \mathrm{InsRay}(z)\}$$

is infinite with positive probability. For the proof of this we need to distinguish when the walk is and is not on $\mathrm{Ray}(z)$. To this end, we introduce a sequence of stopping times defined by $\tau_0 = 0$ and

$$\tau_{n+1} = \inf\{k > \tau_n | X_k \in \mathrm{Ray}(z)\}, \qquad n \geq 0.$$

We will employ two martingale estimates in Lemma 10. Both use Azuma's inequality [1], a version of which we recall in Lemma 9.

LEMMA 9 (Azuma's inequality). *For every $b_1, b_2 \in (0, \infty)$, there exist $b_3, b_4, b_5 \in (0, \infty)$ so that the following holds. If $(Y_n)_{n \geq 0}$ is a sequence of random variables on a probability space with measure $P$ and expectation $E$, and $\tau$ is a (possibly infinite) stopping time w.r.t. that sequence such that $P$-a.s., $Y_0 = 0$, $|Y_{n+1} - Y_n| \leq b_1$ and*

$$(49) \quad E[Y_{n+1} - Y_n | \sigma(Y_m, m \leq n)] \geq b_2 \quad \text{on the event } \{n < \tau\},$$

*then*

$$P[Y_n < b_3 n, n \leq \tau] \leq b_4 e^{-b_5 n} \quad \text{for all } n \geq 0.$$

Lemma 10 will be used in the proof of Lemma 11.

LEMMA 10. *For appropriate $c_{23}, c_{24}$ and $c_{25}$, and all $i \in \{1, 2\}, z \in \partial \mathbb{T}_i$, $x \in \mathrm{InsRay}(z)$ and $n \geq 0$, $\mathbb{P}$-a.s.,*

$$(50) \quad P_{\omega^z}^x[(-1)^{i+1}(X_n - x) \cdot \vec{1} < c_{23} n, T_z \geq n] \leq c_{24} e^{-c_{25} n}$$

*and*

$$(51) \quad P_{\omega^z}^x[T_z = n < \tau_1] \leq c_{24} e^{-c_{25}(u_z(x) \vee n)}.$$

Roughly speaking, (50) says that the speed of the walk is bounded below in the direction $(-1)^{i+1}\vec{1}$, as long as the walk remains in $\mathrm{InsRay}(z)$, whereas (51) bounds the probability of leaving $\mathrm{InsRay}(z) \setminus \mathrm{Ray}(z)$ through the boundary of $\mathrm{InsRay}(z)$.

PROOF OF LEMMA 10. Set $\mathcal{F}_n = \sigma(X_k, k \leq n)$. To demonstrate (50), it suffices to verify that $Y_n = (-1)^{i+1}(X_n - x) \cdot \vec{1}$ satisfies the assumptions of Azuma's inequality, Lemma 9, with $P = P_{\omega^z}^x$ and $\tau = T_z$. Except for (49), all assumptions are obviously satisfied, with $b_1 = 1$. The bound (49) is also satisfied since, on the event $\{T_z > n\}$,

$$(52) \quad \begin{aligned} & E_{\omega^z}^x[(-1)^{i+1}(X_{n+1} - X_n) \cdot \vec{1} | \mathcal{F}_n] \\ & = \sum_{e \in \mathbb{Z}^d, |e|=1} (-1)^{i+1} \omega^z(X_n, X_n + e) e \cdot \vec{1} \\ & \geq (-1)^{i+1} \tfrac{3}{4} r_z(X_n) \cdot \vec{1} - \sum_{e \neq r_z(X_n)} \omega^z(X_n, X_n + e) = \tfrac{3}{4} - \tfrac{1}{4} > 0. \end{aligned}$$

We now demonstrate (51). Since $(X_n)_n$ is a nearest-neighbor walk, $T_z \geq u_z(x)$ and so the statement is trivial for $n < u_z(x)$. Set $Y_n = v_z(X_n) - v_z(x)$. For $n \geq u_z(x)$ on $\{T_z = n < \tau_1\}$, $X_n \notin \mathrm{InsRay}(z)$, and so by (34), $Y_n \leq u_z(X_n) = 0$.

We first consider $x \in \operatorname{InsRay}(z) \setminus \operatorname{Ray}(z)$. For (51) it suffices to check (49) in Azuma's inequality, with $\tau = \tau_1 \wedge T_z$ and $P = P^x_{\omega^z}$, since the other assumptions are obviously satisfied. For this we consider $y \in \operatorname{InsRay}(z) \setminus \operatorname{Ray}(z)$ and estimate the value of $v_z$ at the nearest neighbors of $y$ in terms of $v_z(y)$. For the increment $r_z(y)$,

$$v_z(y + r_z(y)) \stackrel{(33)}{\geq} (n_z(y) + 1)^\beta - |y + r_z(y) - \alpha^{n_z(y)+1}(z)|$$
$$\stackrel{(44)}{\geq} n_z(y)^\beta - |y - \alpha^{n_z(y)}(z)| = v_z(y).$$

Therefore, moving from $y$ to $y + r_z(y)$, which occurs in the environment $\omega^z$ with probability at least $3/4$, does not decrease $v_z$. Similarly,

$$v_z(y + s_z(y)) \stackrel{(33)}{\geq} n_z(y)^\beta - |y + s_z(y) - \alpha^{n_z(y)}(z)|$$
$$\stackrel{(46)}{\geq} n_z(y)^\beta - (|y - \alpha^{n_z(y)}(z)| - 1) = v_z(y) + 1;$$

since a walker at $y$ moves with probability at least $1/5$ to $y + s_z$, $v_z$ increases by 1 with probability at least $1/5$. With probability $1/20$, the walker moves to one of its other neighbors; when doing so, $v_z$ can decrease by at most 1 due to (36). Therefore,

$$(53) \qquad E^x_{\omega^z}[Y_{n+1} - Y_n | \mathcal{F}_n] \geq 1/5 - 1/20 > 0$$

for any $n \geq 0$ on the event $\{n < \tau_1 \wedge T_z\}$, which implies (49) and hence (51) for $x \in \operatorname{InsRay}(z) \setminus \operatorname{Ray}(z)$.

On the other hand, for $x \in \operatorname{Ray}(z)$ and $n = 0$, (51) is trivial. For $x \in \operatorname{Ray}(z)$ and $n > 0$,

$$P^x_{\omega^z}[T_z = n < \tau_1] = \sum_{y \in \mathbb{Z}^d \setminus \operatorname{Ray}(z) : |x-y|=1} \omega^z(x,y) P^y_{\omega^z}[T_z = n - 1 < \tau_1]$$

and we can apply the bound (51) already proved for $y \in \operatorname{InsRay}(z) \setminus \operatorname{Ray}(z)$ to obtain (51) in this case, with a new choice of $c_{24}$. □

Inequality (55) is the main result we will need for the proof of Proposition 12. It follows quickly from (54), which is an analog of (51), but without the restriction on not hitting $\operatorname{Ray}(z)$ before exiting $\operatorname{InsRay}(z)$.

LEMMA 11. *For appropriate $c_{26}$, $c_{27}$, $c_{28}$ and $c_{29}$, and all $i \in \{1,2\}$, $z \in \partial \mathbb{T}_i$, $x \in \operatorname{InsRay}(z)$, and $n \geq 0$,*

$$(54) \qquad P^x_{\omega^z}[T_z = n] \leq c_{26} e^{-c_{27}(u_z(x) \vee n^\beta)}, \qquad \mathbb{P}\text{-a.s.}$$

*and*

$$(55) \qquad E^x_{\omega^z}[T_z; n < T_z < \infty] \leq c_{28} e^{-c_{29}(u_z(x) \vee n^\beta)}, \qquad \mathbb{P}\text{-a.s.}$$

PROOF. By symmetry, we may assume $i = 1$. We will demonstrate (54) by a union bound with four events. Choose $k_n$ so that $\tau_{k_n} \leq n < \tau_{k_n+1}$ and set $\xi_n = (x - z) \cdot \vec{1} + c_{23}n$, where $c_{23}$ is as in Lemma 10. Since $x \in \text{InsRay}(z)$, $\xi_n > 0$. One can then check that

$$P^x_{\omega^z}[T_z = n] \leq \text{I} + \text{II} + \text{III} + \text{IV},$$

where

$$\text{I} = P^x_{\omega^z}[(X_n - x) \cdot \vec{1} < c_{23}n, T_z = n],$$
$$\text{II} = P^x_{\omega^z}[T_z = n < \tau_1],$$
$$\text{III} = P^x_{\omega^z}[\tau_1 < T_z = n, n - \tau_{k_n} \geq \xi_n^\beta/2],$$
$$\text{IV} = P^x_{\omega^z}[(X_n - x) \cdot \vec{1} \geq c_{23}n, \tau_1 < T_z = n, n - \tau_{k_n} < \xi_n^\beta/2].$$

In words, the event in I occurs when, by the time $n$ at which the walk exits $\text{InsRay}(z)$, it has not moved to a much wider part of $\text{InsRay}(z)$. The event in II occurs when, by that time $n$, the walk has not hit $\text{Ray}(z)$. The event in III occurs when, by that time $n$, the walk has hit $\text{Ray}(z)$, but the elapsed time since last visiting $\text{Ray}(z)$ is large. The event in IV occurs when, by that time $n$, the walk has hit $\text{Ray}(z)$, the elapsed time since last visiting $\text{Ray}(z)$ is not large and the walk has moved to a much wider part of $\text{InsRay}(z)$.

We will show that each of these four terms has an upper bound of the form in (54). The bounds for I and II follow directly from Lemma 10, while those for III and IV require some additional work. For I, note that from (50) and $T_z \geq u_z(x)$, it follows that $\text{I} \leq c_{24} e^{-c_{25}(u_z(x) \vee n)}$, which is a stronger bound than required. The same is true for the estimate of II provided by (51). We next bound III. Since $\tau_k$ is the time of the $k$th visit to $\text{Ray}(z)$, $\tau_k \geq k$ for all $k \geq 0$. It follows from this and the Markov property that

$$\begin{aligned}
\text{III} &= \sum_{k=1}^{n} \sum_{l=k}^{\lfloor n - \xi_n^\beta/2 \rfloor} P^x_{\omega^z}[\tau_k = l < T_z = n < \tau_{k+1}] \\
&= \sum_{k=1}^{n} \sum_{l=k}^{\lfloor n - \xi_n^\beta/2 \rfloor} E^x_{\omega^z}[P^{X_l}_{\omega^z}[T_z = n - l < \tau_1]; \tau_k = l < T_z] \\
&\stackrel{(51)}{\leq} \sum_{k=1}^{n} \sum_{l=k}^{\lfloor n - \xi_n^\beta/2 \rfloor} E^x_{\omega^z}[c_{24} e^{-c_{25}(n-l)}; \tau_k = l < T_z] \\
&\leq c_{24} e^{-c_{25}\xi_n^\beta/2} \sum_{k=1}^{n} P^x_{\omega^z}[\tau_k < T_z] \\
&\leq c_{24} n e^{-c_{25}\xi_n^\beta/2}
\end{aligned}$$

(56)

$$\leq c_{26} e^{-c_{30}((x-z)\cdot \vec{1}+n)^\beta}$$

$$\overset{(38)}{\leq} c_{26} e^{-c_{27}(u_z(x) \vee n^\beta)}$$

for appropriate $c_{26}, c_{27}, c_{30}$, where the next to last inequality in (56) follows from the elementary observation that for all $\alpha, \gamma, \delta > 0$ there exists $\eta > 0$ such that for all $s, t \geq 0$,

(57) $$t^\alpha e^{-\gamma(s+t)^\delta} < \eta e^{-\gamma(s+t)^\delta/2}.$$

We demonstrate IV = 0 by showing that the corresponding event cannot occur. We argue by contradiction; on the event in IV,

(58) $$0 = u_z(X_{T_z}) \overset{(34)}{\geq} v_z(X_{T_z}) \overset{(36)}{\geq} v_z(X_{\tau_{k_n}}) - (n - \tau_{k_n}).$$

Since $\tau_1 < n$, we have $k_n \geq 1$ and therefore $X_{\tau_{k_n}} \in \text{Ray}(z)$. Because of (45), $X_{\tau_{k_n}} = \alpha^m(z)$, where $m = n_z(X_{\tau_{k_n}}) = (X_{\tau_{k_n}} - z) \cdot \vec{1}$. So, on the event considered in IV,

$$v_z(X_{\tau_{k_n}})^{1/\beta} = (X_{\tau_{k_n}} - z) \cdot \vec{1}$$
$$= (X_{\tau_{k_n}} - X_n) \cdot \vec{1} + (X_n - x) \cdot \vec{1} + (x - z) \cdot \vec{1}$$
$$\geq \tau_{k_n} - n + \xi_n.$$

Substituting this into (58) and using $n - \tau_{k_n} \leq \xi_n^\beta/2$, we get

$$0 \geq (\xi_n - \xi_n^\beta/2)^\beta - \xi_n^\beta/2 \geq (\xi_n - \xi_n/2)^\beta - \xi_n^\beta/2 = (2^{-\beta} - 2^{-1})\xi_n^\beta > 0,$$

which is a contradiction, so, IV = 0 and we have demonstrated (54).

To obtain (55), note that the left-hand side equals $\sum_{l>n} l P_{\omega^z}^x[T_z = l]$. One can consider separately the cases $u_z(x) \leq n^\beta$ and $u_z(x) > n^\beta$, in the latter case decomposing the sum into $l < u_z^{1/\beta}$ and $l \geq u_z^{1/\beta}$. One can then obtain (55) from (54) and (57) by standard manipulation. $\square$

In the next section we will patch together the different environments $\omega^z$ defined in (47). To do this, it will be useful to introduce some further terminology. Choose $c_{31} > c_{21} \vee 1$ large enough so that for all $n \geq 1$,

(59) $$c_{28} e^{-c_{29} c_{31}^\beta n/c_{21}} \leq \kappa^n,$$

where $c_{21}$ is chosen as in Lemma 8, $c_{28}$ and $c_{29}$ are chosen as in Lemma 11, and $\kappa$ is the ellipticity constant given in (48). For $x \in \mathbb{Z}^d$ and $z \in \partial \mathbb{T}_i$, $i = 1, 2$, define

$$p_z(x) = P_{\omega^z}^x[T_z \leq c_{31} H_i(x)],$$
$$\mathcal{E}_z(x) = E_{\omega^z}^x[T_z; T_z \leq c_{31} H_i(x)],$$
$$\mathcal{E}_z^\infty(x) = E_{\omega^z}^x[T_z; T_z < \infty].$$

Note that

(60)  $\mathcal{E}_z(x) \geq \kappa$ for $x \in \text{InsRay}(z)$ with $d(x, \text{InsRay}(z)^c) = 1$,

since for such $x$, $P^x_{\omega^z}[T_z = 1] \geq \kappa$, and both $c_{31}$ and $H_i(x)$ are at least 1.

Proposition 12 will be used in Section 6. The inequalities in the first line of (61) are elementary; the second line will follow from (55) of Lemma 11.

PROPOSITION 12. *For all $i \in \{1, 2\}, z \in \partial \mathbb{T}_i$ and $x \in \text{InsRay}(z)$,*

(61)
$$\kappa^{u_z(x)} \leq p_z(x) \leq \mathcal{E}_z(x) \leq \mathcal{E}^\infty_z(x),$$
$$\mathcal{E}^\infty_z \leq (c_{28} e^{-c_{29} u_z(x)}) \wedge (\mathcal{E}_z(x) + p_z(x)),$$

$\mathbb{P}$-*a.s., where $c_{28}$ and $c_{29}$ are given in Lemma 10.*

PROOF. Choose a path of length $u_z(x)$ from $x$ to $\text{InsRay}(z)^c$. By (39), $u_z(x) \leq c_{31} H_i(x)$. Therefore, following this path contributes at least probability $\kappa^{u_z(x)}$ to the probability of the event $\{T_z \leq c_{31} H_i(x)\}$, which implies the first inequality of the first line of (61). The other inequalities on this line are immediate from the definition of the quantities involved.

The first part of the inequality in the second line of (61) follows from (55) of Lemma 11, applied to $n = 0$. The second part follows from

$$\mathcal{E}^\infty_z(x) - \mathcal{E}_z(x) = E^x_{\omega^z}[T_z; c_{31} H_i(x) < T_z < \infty] \overset{(55)}{\leq} c_{28} e^{-c_{29}(c_{31} H_i(x))^\beta}$$
$$\overset{(39)}{\leq} c_{28} e^{-c_{29} c_{31}^\beta u_z(x)/c_{21}} \overset{(59)}{\leq} \kappa^{u_z(x)} \leq p_z(x),$$

where the last step employs the first inequality in the first line of (61). $\square$

**6. Patching environments attached to insulated rays.** In this section we prove Theorem 3 by constructing an appropriate random environment $\omega$. The main idea behind the construction of $\omega$ is to choose, for any point $x \in C_i$, among all environments attached to insulated rays covering $x$, the one that "minimizes the probability of exiting the ray." To make this choice locally and thus not destroy the mixing properties inherited by the trees $\mathbb{T}_i$ we have constructed, a slight modification of this idea is actually needed. This is done by minimizing the expectations of exit times from the insulated rays (on the event they are finite).

For $x \in C_1 \cup C_2$, we set $\mathcal{Z}(x) = \{z \in \partial \mathbb{T}_1 \cup \partial \mathbb{T}_2 | x \in \text{InsRay}(z)\}$ and denote by $z(x)$ a leaf $z \in \mathcal{Z}(x)$ that minimizes $\mathcal{E}_z(x)$. (We apply some arbitrary rule, e.g., lexicographic order, to break ties.) Using this, we define, for $x, e \in \mathbb{Z}^d$ with $|e| = 1$,

(62)  $\omega(x, x+e) = \begin{cases} \omega^{z(x)}(x, x+e), & \text{if } x \in C_1 \cup C_2, \\ (2d)^{-1}, & \text{otherwise,} \end{cases}$

where $\omega^z$ is given by (47). Note that $\omega$ inherits the uniform ellipticity of the environments $\omega^z$, with ellipticity constant $\kappa$ given above (48). Moreover, for $x \in C_1 \cup C_2$, we set

$$\mathcal{E}(x) = \mathcal{E}_{z(x)}(x), \qquad \mathcal{E}^\infty(x) = \mathcal{E}^\infty_{z(x)}(x) \quad \text{and} \quad p(x) = p_{z(x)}(x).$$

Because of (60), we will find it useful to employ the stopping time

$$\sigma = \inf\{n \geq 0 | \mathcal{E}(X_n) \geq \kappa\}.$$

The next lemma is the reason for our choice of the environment $\omega$ in (62).

LEMMA 13. *For all $x \in C_1 \cup C_2$, the sequence $(Y_n)_{n \geq 0}$ given by $Y_n = \mathcal{E}(X_{\sigma \wedge n})$ is a supermartingale under $P_\omega^x$ with respect to the filtration $\mathcal{F}_n = \sigma(X_k, k \leq n)$, $n \geq 0$.*

PROOF. Suppose $x \in C_1 \cup C_2$. If $\mathcal{E}(x) \geq \kappa$, then $\sigma = 0$ and the statement is trivial. So, we can assume that $\mathcal{E}(x) < \kappa$. For $y \in C_1 \cup C_2$ with $\mathcal{E}(y) < \kappa$,

$$\begin{aligned}
\mathcal{E}(y) &= \mathcal{E}_{z(y)}(y) \stackrel{(61)}{\geq} \mathcal{E}^\infty_{z(y)}(y) - p_{z(y)}(y) \\
&= E^y_{\omega^{z(y)}}[1 + (T_{z(y)} - 1); T_{z(y)} < \infty] - p_{z(y)}(y) \\
&\geq E^y_{\omega^{z(y)}}[E^{X_1}_{\omega^{z(y)}}[T_{z(y)}; T_{z(y)} < \infty]] = E^y_{\omega^{z(y)}}[\mathcal{E}^\infty_{z(y)}(X_1)] \\
&\stackrel{(61)}{\geq} E^y_{\omega^{z(y)}}[\mathcal{E}_{z(y)}(X_1)].
\end{aligned}$$

Because of $\mathcal{E}(y) < \kappa$ and (60), $d(y, \text{InsRay}(z(y))^c) > 1$. Since the walk is nearest neighbor, this implies that $X_1 \in \text{InsRay}(z(y))$ and, hence, by the definition of $z(X_1)$, $\mathcal{E}_{z(y)}(X_1) \geq \mathcal{E}_{z(X_1)}(X_1)$. Therefore,

$$(63) \qquad \mathcal{E}(y) \geq E^y_{\omega^{z(y)}}[\mathcal{E}_{z(y)}(X_1)] \geq E^y_{\omega^{z(y)}}[\mathcal{E}_{z(X_1)}(X_1)] \stackrel{(62)}{=} E^y_\omega[\mathcal{E}(X_1)].$$

We need to show $E^x_\omega[Y_{n+1}|\mathcal{F}_n] \leq Y_n$. For this, observe that on the event $\{\sigma \leq n\} \in \mathcal{F}_n$, trivially $Y_{n+1} = Y_n$, whereas on the event $\{\sigma > n\}$, by the Markov property,

$$E^x_\omega[Y_{n+1}|\mathcal{F}_n] = E^x_\omega[\mathcal{E}(X_{n+1})|\mathcal{F}_n] = E^{X_n}_\omega[\mathcal{E}(X_1)] \stackrel{(63)}{\leq} \mathcal{E}(X_n) = Y_n.$$

This completes the proof of the lemma. □

We now prove that with positive probability, the random walk $(X_n)_n$ defined by the environment in (62) remains in $C_i$ forever.

PROPOSITION 14. *For $i = 1, 2$, there is $\mathbb{P}$-a.s. some $x \in C_i$, so that*

(64) $$P_\omega^x[X_n \in C_i \text{ for } n \geq 0] > 0.$$

PROOF. For $i = 1, 2$, pick an arbitrary $z \in \partial \mathbb{T}_i$ and set $x = \alpha^n(z) \in \text{Ray}(z)$, where $n$ is large enough so that

(65) $$c_{28} e^{-c_{29} n^\beta} < \kappa$$

for $c_{28}$ and $c_{29}$ chosen as in Lemma 11. By Proposition 12,

(66) $$\mathcal{E}(x) \leq c_{28} e^{-c_{29} u_{z(x)}(x)} \overset{(34)}{\leq} c_{28} e^{-c_{29} n^\beta} \overset{(65)}{<} \kappa.$$

We also require a lower bound on $\mathcal{E}(x)$. Since $(Y_n)_n$ in Lemma 13 is a supermartingale under $P_\omega^x$,

$$\mathcal{E}(x) = E_\omega^x[Y_0] \geq E_\omega^x[Y_n] \geq E_\omega^x[Y_n; \sigma < \infty] = E_\omega^x[\mathcal{E}(X_{\sigma \wedge n}); \sigma < \infty].$$

By Fatou, this implies

$$\mathcal{E}(x) \geq E_\omega^x[\mathcal{E}(X_\sigma); \sigma < \infty] \geq \kappa P_\omega^x[\sigma < \infty].$$

Together with (66), this implies $P_\omega^x[\sigma = \infty] > 0$. On the other hand, by (60), on the event $\{\sigma = \infty\}$, one has $d(X_n, C_i^c) \geq 1$ for all $n$. Therefore, (64) holds. □

We now present the proof of Theorem 3.

PROOF OF THEOREM 3. Define $\omega$ as in (62), with $\beta$ satisfying (24). [Recall that $\beta$ was used throughout the construction of $\omega$, beginning with (25).] By construction, $\omega$ is stationary and is uniformly elliptic, with ellipticity constant at least $\kappa = (20(2d-1))^{-1}$. By Lemma A.5 in the Appendix, one can choose $\beta > 0$ small enough so that $\omega$ is polynomially mixing.

Let $(X_n)_n$ be the RWRE on this environment. We still need to verify that (7) is satisfied. By Proposition 14, with positive probability, $(X_n)_n$ remains forever in $C_i$, $i = 1, 2$, if the RWRE starts at appropriate $x_i \in C_i$. Let $T_{C_i}$ be the exit time of $(X_n)_n$ from $C_i$ (which may be infinite). Since (52) remains valid with the environment $\omega^z$ replaced by $\omega$, a repetition of the proof of (50) shows that

$$P_\omega^{x_i}[(-1)^{i+1}(X_n - x_i) \cdot \vec{1} < c_{23} n, T_{C_i} \geq n] \leq c_{24} e^{-c_{25} n}, \qquad i = 1, 2.$$

By Borel–Cantelli, this implies

$$P_\omega^{x_i}\left[X_n \in C_i \text{ for all } n, \liminf_{n \to \infty} \frac{(-1)^{i+1} X_n \cdot \vec{1}}{n} > \frac{c_{23}}{2}\right] > 0, \qquad i = 1, 2.$$

Since the origin communicates with any $x \in \mathbb{Z}^d$, one obtains (7). □

**7. Open problems.** In this brief section, we mention several open problems. The first involves random forests in $\mathbb{Z}^d$ built from ancestral functions and is motivated by the upper bound in Theorem 2.

OPEN PROBLEM 1. What is the optimal constant $c_1$ in the lower bound (2) of Theorem 1?

There are several natural questions involving RWRE.

OPEN PROBLEM 2. Does the statement of Theorem 3 continue to hold in $d = 2$? If it does, what are the mixing assumptions on the environment that imply the 0–1 law (6)? (Recall that the 0–1 law for i.i.d. environments is proved in [9].)

As mentioned in the Introduction, the following question is still open.

OPEN PROBLEM 3. Prove the 0–1 law (6) for i.i.d. uniformly elliptic environments, when $d \geq 3$.

## APPENDIX

We deferred the demonstration of mixing properties used in the paper to the Appendix; a weaker form of Lemma A.1 was used in the proof of Theorem 2, and Lemma A.5 was used in the proof of Theorem 3. Here, we demonstrate Lemmas A.1 and A.5, and the intermediate Lemmas A.2–A.4 that are used in the proof of Lemma A.5.

We need to extend the notion of polynomial mixing that was introduced in Definition 1, by allowing the set $G$ there to grow with $s$. We will employ the notation $\mathcal{M}_G^b$ introduced in (4).

DEFINITION A.1. Let $\gamma > 0$, $0 \leq \nu < 1$ and $b = (b(x))_{x \in \mathbb{Z}^d}$ be a collection of random variables on a common probability space. Then $b$ is *polynomially $\nu$-mixing (of order $\gamma$)* if, for any fixed $\mu > 0$,

$$\text{(A.1)} \qquad \sup_{s \in \mathbb{Z}^d} \sup_{f \in \mathcal{M}_{\mathcal{B}_s}^b, g \in \mathcal{M}_{\mathcal{B}_s+s}^b} |s|^\gamma |\operatorname{cov}(f,g)| < \infty,$$

where $\mathcal{B}_s = B(0, \mu|s|^\nu)$.

Note that polynomial mixing is the same as polynomial 0-mixing.

Let $a_1 = (a_1(x))_{x \in \mathbb{Z}^d}$ and $a_2 = (a_2(x))_{x \in \mathbb{Z}^d}$ be two independent families of directed ancestral functions that have the same law as the function $a$ constructed in Section 2, where for $a_2$ the direction has been reversed, that is, each $e_j$, $j = 1, \ldots, d$, has been replaced by $-e_j$. The quantities $h_i$ and $H_i$, $i = 1, 2$, are defined analogously to $h$ and $H$ in (1) and (25) by using $a_i$; the quantities $\alpha$ and $\omega$ are given as before by (30) and (62).

We will investigate the mixing properties of the above variables. Our strategy will be to first use i.i.d. random variables $(L_i(x))_{x \in \mathbb{Z}^d, i \in \{1,2\}}$ to construct a realization of the ancestral functions $a_i$; this will allow us to conclude that the pair $(a_1, a_2)$ is polynomially $\nu$-mixing (Lemma A.1). We extend polynomial $\nu$-mixing to the collection $(a_1, h_1, a_2, h_2)$ (Lemma A.2), then to the collections $(a_1, h_1, H_1, a_2, h_2, H_2)$ (Lemma A.3) and $(\alpha, H_1, H_2)$ (Lemma A.4), and finally to the variables $(\omega(x))_{x \in \mathbb{Z}^d}$ (Lemma A.5). In each step, we will use the definitions and appropriate tail estimates to "localize" the random variables that are involved. The details depend on the specific random variables at each step.

The proofs of all five lemmas employ the following elementary inequality: for any measurable functions $f, g, \bar{f}$ and $\bar{g}$ that are bounded in absolute value by 1,

$$
\begin{aligned}
|\operatorname{cov}(f, g)| &\leq |\operatorname{cov}(\bar{f}, \bar{g})| + |\operatorname{cov}(f - \bar{f}, g)| + |\operatorname{cov}(\bar{f}, g - \bar{g})| \\
&\leq |\operatorname{cov}(\bar{f}, \bar{g})| + 4(\mathbb{P}[f \neq \bar{f}] + \mathbb{P}[g \neq \bar{g}]).
\end{aligned}
\tag{A.2}
$$

The first inequality in (A.2) is an immediate consequence of the bilinearity of the covariance function. Throughout this section, in addition to depending on $\beta$ and $d$, all constants $c_i$ are also allowed to depend on $\mu$ and $\nu$.

LEMMA A.1. *For $d \geq 2$ and $0 \leq \nu < 1/d$, the collection $((a_1, a_2)(x))_{x \in \mathbb{Z}^d}$ is polynomially $\nu$-mixing of order $1 - d\nu$.*

PROOF. We may assume that the ancestral functions $a_1$ and $a_2$ have been defined by two independent families $(L_1(x))_x$ and $(L_2(x))_x$ of umbrella lengths.

To prove that $(a_1, a_2)$ is polynomially $\nu$-mixing, we "localize" the variables $a_i$ and show that the localization does not cause damage. Let $\lambda_j^i(x)$ denote the value of $\lambda_j(x)$ that corresponds to the collection $L_i$ as in (12). We set, for $i = 1, 2$, $j = 1, \ldots, d$ and $s, x \in \mathbb{Z}^d$,

$$
\lambda_j^{s,i}(x) = \sup_{y \in B(x, |s|/8) \colon x \in y + (-1)^{i+1} U_{j, L_i(y)}} L_i(y). \tag{A.3}
$$

The random variable $\lambda_j^{s,i}(x)$ is the length of the largest umbrella whose $j$ side passes through $x$ and whose vertex $y$ is contained in $B(x, |s|/8)$. Let $I^{s,i}(x)$, $i = 1, 2$, $s, x \in \mathbb{Z}^d$, be the unique element of $\{1, \ldots, d\}$ for which

$$
\lambda_{I^{s,i}(x)}^{s,i}(x) = \min\{\lambda_j^{s,i}(x) | j = 1, \ldots, d\};
$$

this is the direction with the smallest "protecting" umbrellas. We now set $\bar{a}_i(x) = x + (-1)^{i+1} e_{I^{s,i}(x)}$; this corresponds to the definition of $a(x)$ in Section 2.

Let $(f,g) \in \mathcal{M}_{\mathcal{B}_s}^{(a_1,a_2)} \times \mathcal{M}_{\mathcal{B}_s+s}^{(a_1,a_2)}$ [where the notation is the same as in (4)]. Let $(\bar{f},\bar{g})$ denote the same functions for the collection $(\bar{a}_1,\bar{a}_2)$ instead of $(a_1,a_2)$. To show that (A.1) holds, we use (A.2) to compare $\mathrm{cov}(f,g)$ with $\mathrm{cov}(\bar{f},\bar{g})$.

We will show that $\mathrm{cov}(\bar{f},\bar{g}) = 0$ for $s$ chosen large enough so that $|s| - \mathrm{diam}(\mathcal{B}_s) \geq |s|/2$. To see this, set

$$\hat{\mathcal{B}}_s = \bigcup_{y \in \mathcal{B}_s} B\left(y, \frac{|s|}{8}\right),$$

$$\hat{\mathcal{B}}_s^+ = \bigcup_{y \in \mathcal{B}_s+s} B\left(y, \frac{|s|}{8}\right).$$

Then, for large $s$,

$$d(\hat{\mathcal{B}}_s, \hat{\mathcal{B}}_s^+) \geq |s| - \mathrm{diam}(\mathcal{B}_s) - |s|/4 \geq |s|/4;$$

in particular, $\hat{\mathcal{B}}_s \cap \hat{\mathcal{B}}_s^+ = \varnothing$. Since $\bar{f}$ depends on only those random variables $L_i(x)$ with $x \in \hat{\mathcal{B}}_s$ and since $\bar{g}$ depends on only $L_i(x)$ with $x \in \hat{\mathcal{B}}_s^+$, it follows that $\bar{f}$ and $\bar{g}$ are independent and hence that $\mathrm{cov}(\bar{f},\bar{g}) = 0$.

We next bound $\mathbb{P}[f \neq \bar{f}]$. The functions $f$ and $\bar{f}$ can differ only if there is an $i \in \{1,2\}$, a $j \in \{1,\ldots,d\}$ and an $x \in \mathcal{B}_s$ such that $\lambda_j^{s,i}(x) \neq \lambda_j^i(x)$. In particular, for such $i, j$ and $x$, $\lambda_j^i(x) \geq |s|/8d$, since $\mathrm{diam}(U_{L_i(x)}) \leq dL_i(x)$. Consequently, by Lemma 4, for appropriate $c_{32}$,

$$\mathbb{P}[f \neq \bar{f}] \leq \sum_{i=1}^{2} \sum_{j=1}^{d} \sum_{x \in \mathcal{B}_s} \mathbb{P}\left[\lambda_j^i(x) \geq \frac{|s|}{8d}\right] \leq \frac{16d^2 c_3(\#\mathcal{B}_s)}{|s|} \leq c_{32}|s|^{d\nu-1}.$$

The same bound holds for $\mathbb{P}[g \neq \bar{g}]$. Together with (A.2) and $\mathrm{cov}(\bar{f},\bar{g}) = 0$, this implies (A.1) with $\gamma = 1 - d\nu$ and hence the lemma. $\square$

We extend polynomial $\nu$-mixing from the pair $(a_1,a_2)$ to the collection $(a_1,h_1,a_2,h_2)$.

LEMMA A.2. *For $d \geq 2$ and $0 \leq \nu < 1/d$, the collection $((a_1,h_1,a_2,h_2)(x))_{x \in \mathbb{Z}^d}$ is polynomially $\nu$-mixing of order $(1-d\nu)(d-1)/(2d-1)$.*

PROOF. Fix $\delta = (1 + \nu(d-1))/(2d-1)$ and note that $\delta > \nu$ because $\nu < 1/d$. For any $G \subset \mathbb{Z}^d$ and $s \in \mathbb{Z}^d$, define the event

(A.4) $$A_s(G) = \bigcap_{x \in G} \bigcap_{i=1,2} \{h_i(x) < |s|^{\delta}\}.$$

Also, set $(a,h) = (a_1,h_1,a_2,h_2)$. By stationarity, $\mathbb{P}[A_s(\mathcal{B}_s)] = \mathbb{P}[A_s(\mathcal{B}_s+s)]$, where $\mathcal{B}_s$ is as in Definition A.1. Hence, by (A.2), for any $(f,g) \in \mathcal{M}_{\mathcal{B}_s}^{(a,h)} \times$

$\mathcal{M}^{(a,h)}_{\mathcal{B}_s+s}$,

(A.5)  $\quad |\operatorname{cov}(f,g)| \leq |\operatorname{cov}(f\mathbf{1}_{A_s(\mathcal{B}_s)}, g\mathbf{1}_{A_s(\mathcal{B}_s+s)})| + 4\mathbb{P}[A_s(\mathcal{B}_s)^c].$

To demonstrate that $(a,h)$ is polynomially $\nu$-mixing, we bound the two terms on the right-hand side of (A.5).

To bound the first term, we apply Lemma A.1. For $G \subset \mathbb{Z}^d$, set

$$\mathcal{G}_s(G) = \sigma((a_1(y), a_2(y)), y \in D_s(G)),$$

where

$$D_s(G) = \{y \in \mathbb{Z}^d | d(y, G) \leq |s|^\delta + 1\}.$$

To determine if $h_i(x) < |s|^\delta$ for $x \in G$, it suffices to check whether all branches of $\operatorname{Tree}_i(x)$ terminate within $D_s(G)$. These events are measurable with respect to the ancestral functions $a_i$ restricted to $D_s(G)$ and, hence, are measurable with respect to $\mathcal{G}_s(G)$, so $A_s(G) \in \mathcal{G}_s(G)$. Similarly, the random variables $h_i(x)$, $x \in G$, are determined by $a_i$ restricted to $D_s(G)$. Therefore, setting $G = \mathcal{B}_s$ and $G = \mathcal{B}_s + s$, respectively, it follows that $f\mathbf{1}_{A_s(\mathcal{B}_s)}$ and $g\mathbf{1}_{A_s(\mathcal{B}_s+s)}$ are $\mathcal{G}_s(\mathcal{B}_s)$- and $\mathcal{G}_s(\mathcal{B}_s+s)$-measurable, respectively. It is easy to see that $\operatorname{diam}(D_s(\mathcal{B}_s)) = \operatorname{diam}(D_s(\mathcal{B}_s+s)) \leq c_{33}|s|^\delta$ since $\delta > \nu$, where $c_{33}$ is allowed to depend on $\mu$ but not on $s$, $f$ or $g$. Consequently, by Lemma A.1, for appropriate $c_{34}$,

(A.6)  $\quad |\operatorname{cov}(f\mathbf{1}_{A_s(\mathcal{B}_s)}, g\mathbf{1}_{A_s(\mathcal{B}_s+s)})| \leq c_{34}|s|^{-(1-d\delta)} = c_{34}|s|^{-(d-1)(\delta-\nu)},$

where the last equality follows from the choice of $\delta$.

To bound $\mathbb{P}[A_s(\mathcal{B}_s)^c]$, we note that on $A_s(\mathcal{B}_s)^c$, there exist $x \in \mathcal{B}_s$ and $i \in \{1,2\}$ with $h_i(x) \geq |s|^\delta$. For each $x' \in \operatorname{Ray}(x)$, $h_i(x') \geq h_i(x) \geq |s|^\delta$. Moreover, $\operatorname{Ray}(x)$ must intersect the boundary $\partial \mathcal{B}_s$ of $\mathcal{B}_s$ at some point $y$; thus, on $A_s(\mathcal{B}_s)^c$, $h_i(y) \geq |s|^\delta$. Consequently, for appropriate $c_{35}$,

(A.7)  $\quad \mathbb{P}[A_s(\mathcal{B}_s)^c] \leq 2(\#\partial\mathcal{B}_s)\mathbb{P}[h_1(0) \geq |s|^\delta] \leq c_{35}|s|^{-(d-1)(\delta-\nu)},$

where (3) of Theorem 2 was used in the second inequality. Substitution of (A.6) and (A.7) into (A.5) implies the lemma.  □

We next strengthen the above lemma by including $H_i$. Recall that the definition of $H_i$ in (25) depends on the parameter $\beta \in \mathcal{I}_d = (0, (d-2)/2d)$.

LEMMA A.3. *For $d \geq 3$, fixed $0 < \nu < 1/d$ and all $\beta > 0$ small enough, the collection $((a_1, h_1, H_1, a_2, h_2, H_2)(x))_{x \in \mathbb{Z}^d}$ is polynomially $\nu$-mixing of order $\gamma = (1-d\nu)(d-1)/(2d-1)$.*

PROOF. Assume that $\beta \in \mathcal{I}_d$. We use $(a, h, H)$ as shorthand notation for $(a_1, h_1, H_1, a_2, h_2, H_2)$. Let $(f, g) \in \mathcal{M}^{(a,h,H)}_{\mathcal{B}_s} \times \mathcal{M}^{(a,h,H)}_{\mathcal{B}_s+s}$ for $s \in \mathbb{Z}^d$. For $i = 1, 2$ and $x \in \mathbb{Z}^d$, we set

$$\bar{H}_i(x) = \sup\{h_i(y) | x \in B(y, h_i(y)^\beta), d(y, \mathcal{B}_s) \leq |s|^\nu\},$$
(A.8)
$$\bar{H}_i^+(x) = \sup\{h_i(y) | x \in B(y, h_i(y)^\beta), d(y, \mathcal{B}_s + s) \leq |s|^\nu\}.$$

The quantities $\bar{H}_i$ and $\bar{H}_i^+$ are "localized" versions of $H_i$, which was defined in (25).

Let $\bar{f}$ be defined the same way as $f$, except that one uses the random variables $(a_1, h_1, \bar{H}_1, a_2, h_2, \bar{H}_2)$ instead of $(a, h, H)$. Similarly, let $\bar{g}$ be defined the same way as $g$, except that one uses $(a_1, h_1, \bar{H}_1^+, a_2, h_2, \bar{H}_2^+)$ instead of $(a, h, H)$. Note that $\bar{f}$ (resp. $\bar{g}$) is measurable with respect to the random variables $(a_1, h_1, a_2, h_2)(y)$ with $|y| \leq (\mu + 1)|s|^\nu$ (resp. with $|y - s| \leq (\mu + 1)|s|^\nu$). Therefore, by Lemma A.2,

(A.9) $$\sup_{s, \beta \in \mathcal{I}_d} \sup_{f \in \mathcal{M}^{(a,h,H)}_{\mathcal{B}_s}, g \in \mathcal{M}^{(a,h,H)}_{\mathcal{B}_s+s}} |s|^\gamma |\operatorname{cov}(\bar{f}, \bar{g})| < \infty.$$

To show that (A.1) holds, we use (A.2) to compare $\operatorname{cov}(f, g)$ with $\operatorname{cov}(\bar{f}, \bar{g})$. We still need to bound $\mathbb{P}[f \neq \bar{f}]$ and $\mathbb{P}[g \neq \bar{g}]$. By the definition of $f$ and $\bar{f}$, $f \neq \bar{f}$ can only occur if there exist $i, y$ and $x$ with $i \in \{1, 2\}$, $d(y, \mathcal{B}_s) > |s|^\nu$ and $x \in B(y, h_i(y)^\beta) \cap \mathcal{B}_s$. For such $y$ and $x$, $h_i(y)^\beta \geq |y - x| > |s|^\nu$ and, therefore, $H_i(x) \geq h_i(y) > |s|^{\nu/\beta}$. Consequently,

(A.10)
$$\mathbb{P}[f \neq \bar{f}] \leq \mathbb{P}[H_i(x) \geq |s|^{\nu/\beta} \text{ for some } x \in \mathcal{B}_s \text{ and } i \in \{1, 2\}]$$
$$\leq 2(\#\mathcal{B}_s)\mathbb{P}[H_1(0) \geq |s|^{\nu/\beta}] \overset{(26)}{\leq} c_{36}|s|^{d\nu - ((1-\beta)d-1)\nu/\beta}$$

for appropriate $c_{36}$. [We remind the reader that $c_{36}$ is allowed to depend on $\beta, \mu$ and $\nu$, but not on the specific choice of functions $(f, g) \in \mathcal{M}^{(a,h,H)}_{\mathcal{B}_s} \times \mathcal{M}^{(a,h,H)}_{\mathcal{B}_s+s}$.] For $\beta$ chosen small enough, the right-hand side of (A.10) decays to 0 with exponent larger than $\gamma$. An upper bound for $\mathbb{P}[g \neq \bar{g}]$ is obtained similarly. Together with (A.9) and (A.2), this completes the proof of the lemma. □

The next lemma shows that the triple $(\alpha, H_1, H_2)$ is polynomially $\nu$-mixing. Since the ancestral function $\alpha$ was defined only on $\mathbb{T}_1 \cup \mathbb{T}_2$ [in (30)], we find it convenient to extend the definition, setting $\alpha(x) = \Delta$ for some $\Delta \notin \mathbb{Z}^d$ and all $x \notin \mathbb{T}_1 \cup \mathbb{T}_2$.

LEMMA A.4. *For $d \geq 3$, $\beta > 0$ small enough and all $0 \leq \nu \leq 1/8d$, the collection $((\alpha, H_1, H_2)(x))_{x \in \mathbb{Z}^d}$ is polynomially $\nu$-mixing of order $1/10$.*

PROOF. Since polynomial $\nu$-mixing is monotone in $\nu$, it suffices to show the statement for $\nu = 1/8d$. Let $\beta \in \mathcal{I}_d$ and fix $(f,g) \in \mathcal{M}_{\mathcal{B}_s}^{(\alpha,H_1,H_2)} \times \mathcal{M}_{\mathcal{B}_s+s}^{(\alpha,H_1,H_2)}$ for $s \in \mathbb{Z}^d$. By the definition of $\alpha$, $f$ is a measurable function of the random variables $(a_i(x), H_i(x), \mathbf{1}_{x \in \mathbb{T}_i})_{i=1,2;x \in \mathcal{B}_s}$ and $g$ is a measurable function of $(a_i(x), H_i(x), \mathbf{1}_{x \in \mathbb{T}_i})_{i=1,2;x \in \mathcal{B}_s+s}$.

We proceed to "localize" the variables $\mathbf{1}_{x \in \mathbb{T}_i}$; this will allow us to apply (A.2) the same way as in the previous lemmas. For $(i,j) = (1,2), (2,1)$, set

$$S_i(x) = \bigcap_{0 \leq n \leq |s|^{6\nu}} \bigcap_{y \in B(a_i^n(x), h_i(a_i^n(x))^\beta \wedge |s|^{6\nu})} \{H_j(y) < h_i(a_i^n(x))\}.$$

Let $\bar{f}$ (resp. $\bar{g}$) denote the same function as $f$ (resp. $g$), except that the random variables $\mathbf{1}_{x \in \mathbb{T}_i}$ are replaced by $\mathbf{1}_{S_i(x)}$. Note that $\{x \in \mathbb{T}_i\} \subseteq S_i(x)$, because one recovers the event $\{x \in \mathbb{T}_i\}$ by altering the definition of $S_i(x)$ by removing the restriction $n \leq |s|^{6\nu}$ and not truncating the radius of the ball around $a_i^n(x)$ at $|s|^{6\nu}$. The event $S_i(x)$ is local in the sense that $S_i(x)$ is an element of the $\sigma$-field generated by $(a_i(y), h_i(y), H_j(y))$ with $|y - x| \leq 2|s|^{6\nu}$. (The event $\{x \in \mathbb{T}_i\}$, of course, does not have this property.) Consequently, $\bar{f}$ is measurable with respect to the $\sigma$-field generated by $(a_i(y), h_i(y), H_i(y))$, where $d(y, \mathcal{B}_s) \leq 2|s|^{6\nu}$ and $i \in \{1,2\}$. Similarly, $\bar{g}$ is measurable with respect to the $\sigma$-field generated by $(a_i(y), h_i(y), H_i(y))$, where $d(y, \mathcal{B}_s + s) \leq 2|s|^{6\nu}$ and $i \in \{1,2\}$.

We use the localized functions $\bar{f}$ and $\bar{g}$ together with (A.2) to prove polynomial $\nu$-mixing. For small $\beta$, Lemma A.3 implies that for appropriate $c_{37}$,

$$\begin{aligned}|\operatorname{cov}(\bar{f}, \bar{g})| &\leq c_{37}|s|^{-(1-6d\nu)(d-1)/(2d-1)} \\ &= c_{37}|s|^{-(d-1)/4(2d-1)} \\ &\leq c_{37}|s|^{-1/10},\end{aligned}$$
(A.11)

where $d \geq 3$ is used in the last inequality. To estimate $\mathbb{P}[f \neq \bar{f}]$ and $\mathbb{P}[g \neq \bar{g}]$, we use $\{x \in \mathbb{T}_i\} \subseteq S_i(x)$ and translation invariance to obtain

$$\mathbb{P}[f \neq \bar{f}] \vee \mathbb{P}[g \neq \bar{g}]$$
$$\leq \mathbb{P}\left[\bigcup_{x \in \mathcal{B}_s, i=1,2} S_i(x) \setminus \{x \in \mathbb{T}_i\}\right]$$
$$\leq \sum_{x \in \mathcal{B}_s, i=1,2} \left(\mathbb{P}\left[\bigcup_{n > |s|^{6\nu}} \bigcup_{y \in B(a_i^n(x), h_i(a_i^n(x))^\beta)} \{H_j(y) \geq h_i(a_i^n(x))\}\right]\right.$$
(A.12)
$$\left. + \mathbb{P}\left[\bigcup_{0 \leq n \leq |s|^{6\nu}} \{h_i(a_i^n(x))^\beta > |s|^{6\nu}\}\right]\right)$$

$$\leq (\#\mathcal{B}_s) \sum_{i=1,2} \left( \mathbb{P}\left[ \bigcup_{n>|s|^{6\nu}} \{a_i^n(0) \notin \widetilde{\mathbb{T}}_i\} \right] + \sum_{|y|\leq |s|^{6\nu}} \mathbb{P}[h_i(y) > |s|^{6\nu/\beta}] \right)$$

$$\leq c_{38}(|s|^{d\nu+6\nu((2\beta-1)d+2)} + |s|^{6d\nu+6\nu(1-d)/\beta})$$

for appropriate $c_{38}$. The last inequality uses (28) and (3) of Theorem 2. For $\beta > 0$ small enough, the second term in the right-hand side of (A.12) decays faster than the first. The exponent of $|s|$ in the first term is $-5/8 + 3(1+\beta d)/(2d)$, which is less than $-1/10$ for $d \geq 3$ and $\beta > 0$ small enough. Together with (A.11) and (A.2), this proves the lemma. □

We are finally ready to prove that the environment $\omega$ is polynomially mixing. This result is used in the proof of Theorem 3.

LEMMA A.5. *For $d \geq 3$, with $\beta > 0$ and $\nu > 0$ both small enough, $(\omega(x))_{x \in \mathbb{Z}^d}$ is polynomially $\nu$-mixing of order $1/13$.*

PROOF. Fix $(f,g) \in \mathcal{M}^\omega_{\mathcal{B}_s} \times \mathcal{M}^\omega_{\mathcal{B}_s+s}$. For $G \subseteq Z^d$, denote by $A_s(G)$ the event that $H_i(x) < |s|^{1/8d}$ for all $x \in G$ and $i = 1, 2$. Set $\bar{f} = f\mathbf{1}_{A_s(\mathcal{B}_s)}$ and $\bar{g} = g\mathbf{1}_{A_s(\mathcal{B}_s+s)}$. By Lemma 6,

$$\mathbb{P}[A_s(\mathcal{B}_s)^c] \leq 2(\#\mathcal{B}_s)c_{17}|s|^{(1-(1-\beta)d)/8d} \leq c_{39}|s|^{-(d-1)/8d+d\nu+\beta/8}$$

$$\overset{d\geq 3}{\leq} c_{39}|s|^{-1/12+d\nu+\beta/8} \leq c_{39}|s|^{-1/13}$$

for $\beta > 0$ and $\nu > 0$ small enough, and appropriate $c_{39}$. So, to show polynomial $\nu$-mixing of order $1/13$, it suffices to bound the first term on the right-hand side of (A.2). For this, we will show that $\bar{f}$ and $\bar{g}$ are measurable with respect to $\mathcal{H}_s(\mathcal{B}_s)$ and $\mathcal{H}_s(\mathcal{B}_s + s)$, respectively, where

$$\mathcal{H}_s(G) = \sigma((\alpha(x), H_1(x), H_2(x)), d(x, G) \leq c_{40}|s|^{1/8d})$$

and $c_{40} = 4dc_{31}$, where $c_{31}$ is as in (59). Since the arguments are the same, we will only do this for $\bar{f}$. It will then follow from Lemma A.4 that the first term in (A.2) is bounded above by $c_{41}|s|^{-1/10}$ for appropriate $c_{41}$ not dependent on $f, g$ or $s$.

We note that $A_s(G) \in \mathcal{H}_s(G)$ for $G \subset \mathbb{Z}^d$ and so $A_s(\mathcal{B}_s) \in \mathcal{H}_s(\mathcal{B}_s)$. For $G \subset \mathbb{Z}^d$, write $\mathcal{N}_G$ for the set of functions that are measurable with respect to $\mathcal{H}_s(G)$. Since it is assumed that $f \in \mathcal{M}^\omega_{\mathcal{B}_s}$, to show $\bar{f} \in \mathcal{N}_{\mathcal{B}_s}$, it is clearly enough to show that for $x \in \mathcal{B}_s$,

(A.13) $$\omega(x)\mathbf{1}_{A_s(\mathcal{B}_s)} \in \mathcal{N}_{\mathcal{B}_s}.$$

That is, on the event $A_s(\mathcal{B}_s)$, $\omega(x)$ is a (measurable) function of the random variables that generate $\mathcal{H}_s(\mathcal{B}_s)$.

We first recall how $\omega(x)$ was constructed. Whether $x \in C_1$, $x \in C_2$ or neither holds is determined by $\mathcal{Z}(x)$. [Recall that $z \in \mathcal{Z}(x)$ exactly when $x \in \text{InsRay}(z)$. For $\mathcal{Z}(x) \neq \varnothing$, the direction of $\text{InsRay}(z)$ for any $z \in \mathcal{Z}(x)$ determines whether $x \in C_1$ or $x \in C_2$.] If $\mathcal{Z}(x) = \varnothing$, then the components of $\omega(x)$ all equal $1/2d$. If $\mathcal{Z}(x) \neq \varnothing$, with $x \in C_i$, one computes the random variables $n_z(y), r_z(y)$ and $s_z(y)$ for all $z \in \mathcal{Z}(x)$ and $y \in B(x, c_{31}H_i(x)) \cap \text{InsRay}(z)$, with $c_{31}$ as in (59). From these random variables, one determines the quantities $\omega^z(y)$ as in (47). One then computes $\mathcal{E}_z(x)$, which one uses to determine $z(x)$; one then sets $\omega(x) = \omega^{z(x)}(x)$.

To show (A.13), we therefore proceed as follows.

(a) We show that on $A_s(\mathcal{B}_s)$, for $x \in \mathcal{B}_s$, the random set $\mathcal{Z}(x)$ is a (measurable) function of the random variables that generate $\mathcal{H}_s(\mathcal{B}_s)$, that is, for $z \in \mathbb{Z}^d$, $\mathbf{1}_{z \in \mathcal{Z}(x)} \mathbf{1}_{A_s(\mathcal{B}_s)} \in \mathcal{N}_{\mathcal{B}_s}$.

(b) We next show that on $A_s(\mathcal{B}_s)$, $x \in \mathcal{B}_s \cap C_i$, $z \in \mathcal{Z}(x)$ and $y \in B(x, c_{31}H_i(x)) \cap \text{InsRay}(z)$, the random variables $n_z(y), r_z(y)$ and $s_z(y)$ are functions of the random variables that generate $\mathcal{H}_s(\mathcal{B}_s)$.

(c) Finally, we show that on $A_s(\mathcal{B}_s)$, $x \in \mathcal{B}_s \cap C_i$ and $z \in \mathcal{Z}(x)$, $\mathcal{E}_z(x)$ is a function of the random variables that generate $\mathcal{H}_s(\mathcal{B}_s)$.

The following inclusion, whose justification we defer to the end of the proof, is used for all three steps. For all $x \in \mathcal{B}_s$, $z \in \mathcal{Z}(x)$ and $y \in B(x, c_{31}|s|^{1/8d}) \cap \text{InsRay}(z)$,

$$(\text{A.14}) \quad \{\alpha^n(z) | 0 \le n \le n_z(y)\} \subseteq B(x, c_{40}|s|^{1/8d}/2) \quad \text{on } A_s(\mathcal{B}_s).$$

In particular, (a) is an immediate consequence of (A.14) with $y = x$ and the definition of $\text{InsRay}(z)$.

To see (b), first note that by (A.14), on $A_s(\mathcal{B}_s)$, the variables $n_z(y)$, $\alpha^{n_z(y)}(z)$ and $\alpha^{n_z(y)+1}(z)$ are functions of the random variables that generate $\mathcal{H}_s(\mathcal{B}_s)$. [The set $B(x, c_{40}|s|^{1/8d}/2)$ was enlarged to $B(x, c_{40}|s|^{1/8d})$ to include $\alpha^{n_z(y)+1}$ in $\mathcal{H}_s(\mathcal{B}_s)$.] Since $r_z(y)$ and $s_z(y)$ are determined by $\alpha^{n_z(y)}(z)$ and $\alpha^{n_z(y)+1}(z)$, (b) follows.

To see (c), recall that for $x \in C_i$, $\mathcal{E}_z(x) = E_{\omega^z}^x[T_z; T_z \le c_{31}H_i(x)]$. The RWRE is nearest neighbor and so, starting at $x$, will not escape $B(x, c_{31}|s|^{1/8d})$ by time $c_{31}|s|^{1/8d}$. On $A_s(\mathcal{B}_s)$, $H_i(x) \le |s|^{1/8d}$. Consequently, on $A_s(\mathcal{B}_s)$, $\mathcal{E}_z(x)$ is a function of $\omega^z$ on $B(x, c_{31}|s|^{1/8d})$ and of $\text{InsRay}(z)$. By (A.14), (47) and (b), the claim in (c) holds.

It only remains to show (A.14). First observe that on $A_s(\mathcal{B}_s)$, for any $z \in \mathcal{Z}(x)$,

$$(\text{A.15}) \quad \begin{aligned} |z - x| &\le |z - \alpha^{n_z(x)}(z)| + |\alpha^{n_z(x)}(z) - x| \le n_z(x) + n_z(x)^\beta \\ &\le 2n_z(x) \le 2(H_1(x) \vee H_2(x)) \le 2|s|^{1/8d}. \end{aligned}$$

The second inequality follows from the definitions of $n_z(x)$ and $\mathrm{InsRay}(z)$ in (35) and (32); the fourth inequality follows from the definition of $H$ in (25) and the inclusion $x \in B(\alpha^{n_z(x)}(z), n_z(x)^\beta)$.

Since $B_\infty(x, 2c_{31}|s|^{1/8d}) \subseteq B(x, c_{40}|s|^{1/8d}/2)$ and since the path $\alpha^n(z)$, $n = 0, \ldots, n_z(y)$, is directed, it suffices to show that both endpoints $z$ and $\alpha^{n_z(y)}(z)$ are contained in $B_\infty(x, 2c_{31}|s|^{1/8d})$. This holds for $z$ because of (A.15). For $\alpha^{n_z(y)}(z)$, first note that

$$n_z(y)^\beta \geq |\alpha^{n_z(y)}(z) - y|$$
$$\geq |\alpha^{n_z(y)}(z) - z| - |z - x| - |x - y|$$
$$\geq n_z(y) - (2 + c_{31})|s|^{1/8d}.$$

The first inequality follows from the definition of $\mathrm{InsRay}(z)$; the last inequality follows from (A.15), since $y \in B(x, c_{31}|s|^{1/8d})$. Since $\beta < 1$, this implies $n_z(y)^\beta \leq |s|^{1/8d}$ for $|s|$ large. Therefore,

$$|\alpha^{n_z(y)}(z) - x|_\infty \leq |\alpha^{n_z(y)}(z) - x| \leq |\alpha^{n_z(y)}(z) - y| + |y - x| \leq 2c_{31}|s|^{1/8d}$$

for $|s|$ large, where we used $y \in \mathrm{InsRay}(z)$ in the last inequality. This demonstrates (A.14) and completes the proof of the lemma. $\square$

REMARK. The only place where the explicit structure of the function $a$ is used is in the proof of Lemma A.1. Given any ancestral functions $\bar{a}_i$, $i = 1, 2$, that are directed in opposite directions and for which the conclusion of Lemma A.1 holds, Lemmas A.2–A.5 will continue to hold. In particular, the environment $\bar{\omega}$ that is constructed from such $\bar{a}_1, \bar{a}_2$ by using the pruning and insulation recipe leading up to (62) will be polynomially $\nu$-mixing of order $1/13$.

M. BRAMSON
SCHOOL OF MATHEMATICS
UNIVERSITY OF MINNESOTA
MINNEAPOLIS, MINNESOTA 55455
USA
E-MAIL: bramson@math.umn.edu

O. ZEITOUNI
SCHOOL OF MATHEMATICS
UNIVERSITY OF MINNESOTA
MINNEAPOLIS, MINNESOTA 55455
USA
AND
DEPARTMENT OF MATHEMATICS
   AND DEPARTMENT OF ELECTRICAL ENGINEERING
TECHNION
HAIFA 32000
ISRAEL
E-MAIL: zeitouni@math.umn.edu
URL: www-ee.technion.ac.il/~zeitouni

M. ZERNER
MATHEMATISCHES INSTITUT
UNIVERSITÄT TÜBINGEN
AUF DER MORGENSTELLE 10
72076 TÜBINGEN
GERMANY
E-MAIL: martin.zerner@uni-tuebingen.de
URL: www.mathematik.uni-tuebingen.de/~zerner